
\documentstyle{amsppt} \magnification1200 \font\scr=eusm10
\loadbold \define\sctheta{\boldsymbol\varTheta}
\define\th{^{\text{th}}} \define\script#1{\text{\scr #1}}
\def\d#1{\tfrac{\partial}{\partial#1}}
\define\comp{\lower2.75pt\hbox{$^\circ$}}
\define\sigmabar{\sigma\kern-6pt\lower1pt\hbox{\vrule depth.2pt
width4pt height.2pt}} \NoBlackBoxes

\voffset=-0.5in

\phantom{}\vskip-\baselineskip \topmatter \title Siegel modular
forms generated by\\invariants of cubic hypersurfaces\endtitle
\rightheadtext{Siegel Modular Forms} \author C. McCrory, T.
Shifrin and R. Varley\endauthor \affil University of Georgia
\endaffil \address University of Georgia, Athens Georgia 30602
\endaddress \email clint\@joe.math.uga.edu,
shifrin\@joe.math.uga.edu, robert\@joe.math.uga.edu \endemail
\thanks Partial support provided by NSF grants DMS-9106938,
DMS-9208282 and DMS-9305857\endthanks \endtopmatter

We give a geometric derivation of Schottky's equation in genus
four for the period matrices of Riemann surfaces among all period
matrices.  The equation arises naturally from the singularity
theory of the Gauss map on the theta  divisor, and thus
generalizes for any genus $g \ge 4$ to a certain ideal of Siegel
modular forms vanishing on period matrices of Riemann surfaces.
This ideal is generated by modular forms associated to the
invariants of cubic forms in $g - 1$ variables which vanish on
the Fermat cubic.

For each $g$-dimensional principally polarized abelian variety
$(A,\Theta )$ together with a marked odd theta characteristic
$\xi$, consider the theta function $\vartheta[\xi](z, \Omega)$
and its  Taylor expansion with respect to $z$ about the origin,
$$ \vartheta[\xi](z, \Omega) = \ell(z, \Omega) + m(z, \Omega) +
\dots, $$ where  $\ell = \ell_{\xi}$ and $m = m_{\xi}$ are,
respectively, linear and cubic forms in $z$.  Classical
transformation formulas show that all of the ideals $(\ell)$,
$(\ell, m)$, $\dots$, transform by the same automorphy factor
under the action on $\Omega$ of a suitable finite index subgroup
of $\text{Sp}(2g, \Bbb Z)$.  The main idea---suggested by the
singularity theory of the Gauss map---is simply to use the ideal
$(\ell, m) \subset \Bbb C[z]$ as an invariant of $(A, \Theta,
\xi)$.  When the linear form $\ell$ is not identically zero, we
can consider the restriction $\overline{m}$ of the cubic form $m$
to the hyperplane $\ell(z) = 0$. This cubic form $\overline{m}$
is the key object studied in this paper. Our main observation is
that if $(A, \Theta)$ is the Jacobian of a generic curve of genus
$g$, then for any choice of odd theta characteristic,
$\overline{m}(z)$ is a Fermat cubic, the sum of $g - 1$ cubes.

Geometrically, the marking of an odd theta characteristic $\xi$
of  a principally polarized abelian variety $(A, \Theta)$
determines a symmetric theta divisor $\Theta[\xi]$ on which the
origin of $A$ has odd multiplicity.  When $\Theta[\xi]$ is smooth
at $0$, its local tangent hyperplane section at $0$ is singular
and of odd multiplicity. In particular, if the cubic form
$\overline{m}(z)$ is not identically zero, then it defines the
projectivized tangent cone  $V$ of this singularity. The
projectivized tangent space $\Bbb E = \Bbb PT_0(\Theta[\xi])$ is
a hyperplane in the canonical projective space $\Bbb P^{g-1} =
\Bbb PT_0(A)$.  Thus we associate to $(A, \Theta, \xi)$ the cubic
hypersurface $V \subset \Bbb E \cong \Bbb P^{g-2}$, and for a
generic Jacobian, $V$ is projectively equivalent to the standard
Fermat cubic.

The main result of this paper is a precise way to convert
invariants vanishing on the Fermat cubic into modular forms
vanishing on Jacobians. If $\varphi$ is any homogeneous
polynomial invariant of cubic forms in $g - 1$ variables, then
there is a natural globalization $G_{\varphi}$ so that the
equation $\varphi(\overline{m}) = 0$ is expressed by the
vanishing of a Siegel modular form $G_{\varphi}(\overline{m})$.
Therefore, if the invariant $\varphi$ vanishes on the Fermat
cubic, then $G_{\varphi}(\overline{m})$ vanishes on Jacobians.
The {\it Fermat ideal\/} is the ideal generated by the Siegel
modular forms resulting from all invariants vanishing on the
Fermat cubic.

We prove in genus $g = 4$ that the zero set of the Fermat ideal
(which is principal in this case) is exactly the closure of the
locus of period matrices of genus four Riemann surfaces, and we
deduce that the generator is Schottky's equation (up to a
constant multiple). (References for genus four Schottky are
\cite{I4}, \cite{Fr1}, \cite{E}.)  The Fermat ideal consists
entirely of cusp forms, and hence for $g \ge 5$ it cannot
coincide with the ideal of Schottky-Jung equations.  (For
Schottky-Jung, cf\. \cite{F1}, \cite{F2}, \cite{vG}, \cite{Do3},
\cite{C}, \cite{M2}, \cite{R-F}.)  In general we obtain two loci
in $\Cal A_g$, the level $1$ moduli space for principally
polarized abelian varieties, corresponding to  the vanishing of
the Fermat ideal for {\it some\/} odd theta characteristic or for
{\it all\/} odd theta characteristics.  These loci are algebraic
subsets of $\Cal A_g$ containing the Jacobians; we call them the
{\it big\/} and {\it small\/} Fermat loci, by analogy with the
Schottky-Jung theory \cite{Do1}.

It is well known that interesting invariants of a principally
polarized abelian variety $(A, \Theta)$ can be obtained from the
common zeros of various terms in the Taylor expansions in $z$ of
theta functions; cf\. \cite{ACGH, Cor\. p\. 232}, \cite{K1}.  As
far as we know, the  particular combination $\ell, m$ (defined by
an odd theta characteristic $\xi$) used here has not been
investigated before (apart from our previous work \cite{AMSV2}).
In fact, the condition that $\overline{m}$ be a Fermat cubic can
be viewed as a necessary algebraic condition for a theta divisor
to be of translation type (cf\. Prop\. (2.3.4), \cite{B, p\. 105
(c)}, \cite{I5, p\. 167}, \cite{L}, \cite{M2}).

In \S 1 we define the various moduli spaces of abelian varieties
that we will use, and we recall the transformation laws for the
theta function $\vartheta[\xi](z, \Omega)$.  In particular we
recall the space $\Cal A_g^{(4, 8)}$, its fundamental positive
line bundle $P$ (which is well known in the theory of theta
constants) and the dual line bundle $N$, with $P^2 \cong \Lambda
= \det H$, where  $H$ is the Hodge vector bundle. The
transformation law for $\vartheta[\xi]$ shows that the automorphy
factor for the line bundle $P$ governs every term, modulo all the
previous terms, in the Taylor expansion of $\vartheta[\xi]$ with
respect to $z$.  Thus, for each odd theta characteristic $\xi$,
the linear term of $\vartheta[\xi]$ defines a homomorphism
$\ell_{\xi}: N \rightarrow H$.

In \S 2 we introduce an open set  $\Cal U_\xi$ over which the
quotient sheaf $H/\ell_{\xi}(N)$ is a rank $g - 1$ vector bundle
$\Cal E_{\xi}^*$, so that the cubic term of $\vartheta[\xi]$
defines a homomorphism $\overline{m}_{\xi}: N \rightarrow S^3\Cal
E_{\xi}^*$ to a vector bundle of cubic forms in $g - 1$
variables, and we check that the complement of $\Cal U_\xi$ is
Zariski-closed of codimension at least $2$.  We observe that, for
every odd theta characteristic $\xi$ of the Jacobian of a generic
curve of genus $g$, the cubic form $\overline{m}_{\xi}(z)$ is a
Fermat cubic (extending the genus three and four cases treated in
\cite{MSV1}, \cite{AMSV2}).

In \S 3 we show that a degree  $d$  polynomial invariant
$\varphi$  of cubic forms in $g - 1$ variables yields, for every
rank $g - 1$ vector bundle $\Cal E$, a degree $d$ polynomial
mapping $G_{\varphi}: S^3\Cal E^* \rightarrow (\Lambda^{g-1}\Cal
E^*)^w$, where $3d = (g - 1)w$. In \S 4 we conclude that each
degree $d$ polynomial mapping $G_{\varphi} \comp
\overline{m}_{\xi}: N \rightarrow (\Lambda^{g-1} \Cal
E_{\xi}^*)^w$ can be represented (in a canonical way) by a Siegel
modular form of a certain weight for the group $\Gamma(4, 8)$.
Combining this general construction with the property that  for
Jacobians the cubic form $\overline{m}_{\xi}(z)$ is a Fermat
cubic, we obtain a simple way to construct Siegel modular forms
vanishing on period matrices of Jacobians. We recover the
Schottky equation in genus $g=4$ and we show that these modular
forms also vanish on products for every genus $g \ge 4$.  In \S 5
we indicate some open questions.

We would like to thank A.~Beauville, O.~Debarre, E.~Izadi, T.~
Johnsen, R.~ Salvati~Manni and R.~Smith for helpful comments.
The first and third authors are grateful for hospitality at MSRI,
where the main results of this paper were first presented in
December 1992.

\head 1. Theta transformation laws \endhead

\subhead 1.1 \quad The moduli space $\Cal A_g(4,8)$\endsubhead
\vskip.5\baselineskip

Let $\Cal A_g$ be the moduli space of $g$-dimensional principally
polarized abelian varieties $(A, \Theta)$.  Here the theta
divisor $\Theta$ is specified only up to translation in $A$; a
symmetric representative of $\Theta$ is specified only up to
translation by a point of order $2$.  The moduli space of primary
interest in this paper is $\tilde{\Cal A}_g$, the isomorphism
classes of $(A, \Theta,\xi)$, $g$-dimensional principally
polarized abelian varieties together with an odd theta
characteristic $\xi$. Recall that the $2^{2g}$ theta
characteristics of a principally polarized abelian variety
correspond precisely to the symmetric representatives of the
theta divisor, and that the parity of a characteristic $\xi$ is
the same as the parity of $\text{mult}_0(\Theta[\xi])$, the
multiplicity of the corresponding symmetric theta divisor at the
origin of the abelian variety.  Thus, a point $(A, \Theta, \xi)$
of $\tilde{\Cal A}_g$ is simply a $g$-dimensional abelian variety
$A$ with a distinguished symmetric theta divisor $\Theta[\xi]$ on
which the origin has odd multiplicity.

Faced with the usual technical problems (related to coarse moduli
and fixed points for the action of the modular group
$\text{Sp}(2g, \Bbb Z)$ on the Siegel upper half space $\Cal
H_g$), we will do most of our calculations on a finite cover
$\Cal A'_g \rightarrow \tilde{\Cal A}_g$, where $\Cal A'_g = \Cal
H_g/\Gamma'$ for a suitable subgroup $\Gamma'$ of $\text{Sp}(2g,
\Bbb Z)$.  In fact, among the spaces that will serve our
purposes, there is a distinguished one in the theory of theta
transformation laws.  This particularly useful space is $\Cal
A_g^{(4, 8)} = \Cal H_g/\Gamma(4, 8)$, introduced by Igusa
\cite{I1}.  We abbreviate $\Cal A_g^{(4, 8)}$ to $\Cal A'_g$ and
use $(A, \Theta, \xi')$ to denote a point of $\Cal A'_g$ with
image $(A, \Theta, \xi)$ in $\tilde{\Cal A}_g$.  We will not need
the precise moduli interpretation of $\xi'$, a level $(4, 8)$
structure (cf\. \cite{M-F, app\. B to Ch\. 7, pp\. 191-199, esp\.
p\. 195}); it would be sufficient to work on the level $8$ moduli
space $\Cal A_g^{(8)} = \Cal H_g/\Gamma(8)$, or with any
principal congruence level  $n$ such that $8|n$. We will use
primarily the following properties: $\Cal A'_g$ is smooth, maps
finitely onto $\tilde{\Cal A}_g$, and carries a family of abelian
varieties together with a family of theta divisors corresponding
to the marked odd theta characteristic $\xi$. (It is convenient
that also the transformation law for the theta function
$\vartheta[\xi]$ becomes quite simple, as we will see below.)

In summary, the relevant moduli spaces and maps are: $$ \Cal H_g
\rightarrow \Cal A'_g \rightarrow \tilde{\Cal A}_g \rightarrow
\Cal A_g. $$ We must choose a map $\pi_\xi: \Cal A'_g \rightarrow
\tilde{\Cal A}_g$ factoring the Galois covering $\Cal A'_g
\rightarrow \Cal A_g$ through the canonical (but non-Galois)
covering $\tilde{\Cal A}_g \rightarrow \Cal A_g$.  The level
structure $\xi'$ which is specified on each element  $(A, \Theta,
\xi') \in \Cal A'_g$ determines a ``framing''  of the set of
theta characteristics of $(A, \Theta)$, so it makes sense to
speak of fixing (universally) an odd theta characteristic $\xi$
on every $(A, \Theta, \xi') \in \Cal A'_g$, and such a choice
defines a map $\pi_\xi$ as above.  To be more precise, we first
pass from $\Cal A'_g$ to the level $2$ moduli space $\Cal
A_g^{(2)}$, which is Galois over $\Cal A_g$ with group
$\text{Sp}(2g, \Bbb Z)/\Gamma(2) \cong \text{Sp}(2g, \Bbb Z/2)$.
A level $2$ structure of a principally polarized abelian variety
$(A, \Theta)$ corresponds to a symplectic isomorphism $(\Bbb
Z/2)^{2g} \rightarrow H_1(A, \Bbb Z/2)$, where $(\Bbb Z/2)^{2g}$
has the standard $\Bbb Z/2$-symplectic form, and $H_1(A, \Bbb
Z/2)$ has the mod $2$ reduction of the alternating Riemann form
on $H_1(A, \Bbb Z)$ determined by the polarization.  Thus a level
$2$ structure allows the specification of an odd theta
characteristic of $(A, \Theta)$ simply by a suitable $\Bbb
Z/2$-quadratic form on $(\Bbb Z/2)^{2g}$ (cf\. \cite{M1, \S 2},
\cite{I3, pp\. 209-210}).

\vskip.5\baselineskip \subhead 1.2\quad Automorphy factors
\endsubhead \vskip.5\baselineskip

We recall the basic theory of automorphy factors (cf\. \cite{G,
\S 5.1, esp\. (5.1.1) p\. 198, p\. 199}, \cite{L-B, App\. B, pp\.
412-414}).  Suppose that a group $\Gamma$ acts properly
discontinuously (as in \cite{G, (D) p\. 203}, \cite{M3, (c), p\.
180}) on a topological space $X$ by a map $\Gamma \times X
\rightarrow X$, $(\gamma, x) \mapsto \gamma \cdot x$, and let $p:
X \rightarrow Y = X/\Gamma$ be the quotient map. Let $E$ be a
rank $r$ complex vector bundle on $Y$ and $\tilde{E} = p^*(E)$
the pullback to $X$. Then $\Gamma$ acts in a natural way on the
vector bundle $\tilde{E}$ and its space of sections $H^0(X,
\tilde{E})$, and $H^0(Y, E) \overset\sim\to\rightarrow (H^0(X,
\tilde{E}))^{\Gamma}$; i.e., sections of $E$ over $Y$ correspond
to invariant sections of $\tilde{E}$ over $X$.  Now suppose that
$\tilde{E}$ is trivial on $X$ and let $\varphi: \tilde{E}
\rightarrow X \times \Bbb C^r$ be a trivialization. Then via the
trivialization there is an action of $\Gamma$ on $X \times  \Bbb
C^r$, $\gamma \cdot (x, w) = \varphi(\gamma \cdot(\varphi^{-1}(x,
w))) = (x', w')$, where  $x' = \gamma \cdot x$ and $w' =
\mu(\gamma, x) w$. The mapping $$ \mu: \Gamma \times X
\rightarrow \text{GL}_r(\Bbb C) $$ is called the {\it automorphy
factor\/} for $E$ (associated to  the trivialization $\varphi$ of
$\tilde{E}$ on $X$) and satisfies the cocycle condition $$
\mu(\gamma'\gamma, x) = \mu(\gamma', \gamma \cdot x) \mu(\gamma,
x) $$ for all $\gamma, \gamma' \in \Gamma$, $x \in X$.  Under the
trivialization $\varphi$, a section $s$ of $\tilde{E}$ over $X$
corresponds to a $\Bbb C^r$-valued function $f$ on $X$ by the
formula $$ \varphi(s(x)) = (x, f(x)). $$ The action of $\Gamma$
on sections of $\tilde{E}$ can then be written $$ (\gamma \cdot
f)(x) = \mu(\gamma, \gamma^{-1} \cdot x) f(\gamma^{-1}\cdot x).
$$ It follows that $f$ is invariant under this action of $\Gamma$
if and only if for every $\gamma \in \Gamma$ and $x \in X,
f(\gamma \cdot x) = \mu(\gamma, x)f(x)$.

Thus if $E$ is a rank $r$ complex vector bundle on $Y =
X/\Gamma$, then a trivialization of the pullback $\tilde{E}$ on
$X$ induces a canonical isomorphism $$ H^0(Y, E)
\overset\sim\to\rightarrow \{f: X \rightarrow \Bbb C^r\ |\
f(\gamma \cdot x) = \mu(\gamma, x)f(x) \text{ for all } \gamma
\in \Gamma, x \in X\}, $$ where  $\mu: \Gamma \times X
\rightarrow \text{GL}_r(\Bbb C)$ is the automorphy factor for
$E$.  By the same reasoning for each open subset $V$ of $Y$, we
get a description of the sheaf of sections of $E$ over $Y$.

We will work with various vector bundles on $\Cal A'_g$; their
natural descriptions are by automorphy factors over $\Cal H_g$,
so we can carry out the constructions in terms of sheaves, i.e.,
sections of the bundles over variable open subsets of $\Cal
A'_g$.  Thus we work mostly in the holomorphic category, with the
analytic structure sheaf $\Cal O$ and analytic $\Cal O$-modules.

\vskip.5\baselineskip \subhead 1.3\quad The Hodge vector bundle
\endsubhead \vskip.5\baselineskip

We describe the modular group action on Siegel upper half-space
$\Cal H_g$ and on the product $\Bbb C^g \times \Cal H_g$. The
space $\Cal H_g$ of all $g \times g$ {\it period matrices\/} is
the set of all $g \times g$ complex matrices $\Omega$ such that
$^t\Omega = \Omega$ and $\text{Im}(\Omega) > 0$.  For $\gamma \in
\text{Sp}(2g,\Bbb Z)$, we write $$ \spreadlines{-.03in} \gamma =
\pmatrix A & B \\ C & D \endpmatrix $$ with  $g \times g$ blocks
$A$, $B$, $C$ and $D$ satisfying $$ \spreadlines{-.03in}
{\vphantom{\pmatrix A\\C\endpmatrix}}^t\!\pmatrix A & B\\ C & D
\endpmatrix \pmatrix 0 & I \\ -I & 0 \endpmatrix \pmatrix A & B
\\ C & D \endpmatrix  = \pmatrix 0 & I \\ -I & 0 \endpmatrix . $$
Then for $\Omega \in \Cal H_g$, let $$ \spreadlines{-.03in}
\gamma \cdot \Omega = (A\Omega + B)(C \Omega + D)^{-1},$$ and for
$(z, \Omega) \in \Bbb C^g \times \Cal H_g$, let $$
\spreadlines{-.03in} \gamma \cdot (z, \Omega) = (^t(C\Omega +
D)^{-1}z, (A\Omega + B) (C\Omega + D)^{-1}) $$ (cf\. \cite{I2,
p\. 226}, \cite{M2, pp\. 71-72}, \cite{M3, p\. 177}, \cite{R-F,
pp\. 86-87}).

For the moduli-theoretic meaning of these actions of
$\text{Sp}(2g, \Bbb Z)$, see \cite{M3, pp\. 171-177, 184-185};
$\Cal H_g$ can be regarded as a fine moduli space for principally
polarized abelian varieties $(A, \Theta)$ with ``level $\infty$
structure," i.e., a symplectic isomorphism $\alpha: \Bbb Z^{2g}
\rightarrow H_1(A, \Bbb Z)$.  Then $\gamma \in \text{Sp}(2g, \Bbb
Z) \subset \text{GL}_{2g}(\Bbb Z)$ acts by $\gamma \cdot(A,
\Theta, \alpha) = (A, \Theta, \alpha \comp \gamma^{-1})$,  and
one obtains the standard formula above for $\gamma \cdot \Omega$
after an automorphism of $\text{Sp}(2g, \Bbb Z)$ \cite{M3, p\.
174}.

Over $\Cal H_g$ there is the universal family of abelian
varieties $\{A_{\Omega} = \Bbb C^g/(\Bbb Z^g + \Omega \Bbb
Z^g)\}$, and the action of $\text{Sp}(2g, \Bbb Z)$ on $\Cal H_g$
lifts canonically to the total space of this family. The family
of abelian varieties defines a vector bundle $\tilde{\Cal V}
\rightarrow \Cal H_g$ of tangent spaces $\{T_0A_{\Omega}\}$, with
an induced action of $\text{Sp}(2g, \Bbb Z)$.  This vector bundle
is trivial; i.e., $\tilde{\Cal V} \cong \Bbb C^{g} \times \Cal
H_g$ over $\Cal H_g$, and the action of $\text{Sp}(2g, \Bbb Z)$
on the bundle corresponds to the automorphy factor $\mu(\gamma,
\Omega) = {}^t(C\Omega + D)^{-1}$; that is, $\gamma \cdot (z,
\Omega) = (\mu(\gamma, \Omega)z, \gamma \cdot \Omega)$, for
$\gamma \in \text{Sp}(2g, \Bbb Z)$ and $(z, \Omega) \in \Bbb C^g
\times \Cal H_g$.

Let $\Cal V = (\Bbb C^g \times \Cal H_g)/\Gamma'$, a rank $g$
vector bundle on $\Cal A'_g = \Cal H_g/\Gamma'$.  The fibre of
$\Cal V$ at $y = (A, \Theta, \xi') \in \Cal A'_g$ is $V_y =
T_0A$. Thus $V_y^* = T_0^*A$, so $\Cal V = H^*$, the dual of the
{\it Hodge vector bundle\/} $H$ on $\Cal A'_g$.  For $\gamma =
\left(\smallmatrix A & B \\ C & D \endsmallmatrix\right) \in
\text{Sp}(2g, \Bbb Z)$ and $\Omega \in \Cal H_g$, let
$\eta(\gamma, \Omega) = C\Omega + D$.  Then $\eta$ is the
automorphy factor for $H$ and $\det\eta(\gamma, \Omega) =
\det(C\Omega + D)$ is the automorphy factor for the {\it Hodge
line bundle\/} $\Lambda = \det H$.

\vskip.5\baselineskip \subhead 1.4\quad Transformation laws
\endsubhead \vskip.5\baselineskip

We consider transformation laws for theta functions
$\vartheta[\xi](z, \Omega)$ and regard the relevant automorphy
factor over $\Cal H_g$ as the description of a fundamental line
bundle on the moduli space $\Cal A'_g = \Cal A_g^{(4, 8)}$. First
we recall the classical transformation theory for theta functions
$\vartheta[\xi](z, \Omega)$ with characteristics.  The
characteristic $[\xi]$ is specified by two $g$-entry column
vectors  $a, b$ of zeros and ones.  As a function of $z \in \Bbb
C^g$ for fixed $\Omega \in \Cal H_g$ and $\xi \in \{0, 1\}^{2g}$,
$$ \vartheta[\xi](z, \Omega) = \sum_{n \in \Bbb Z^g}
\exp\left(\pi i\big({}^t(n + \tfrac a2) \Omega(n + \tfrac a2) +
2\ ^t(n + \tfrac a2)(z + \tfrac b2)\big)\right), $$ where
$\exp(w) = e^w$ for $w \in \Bbb C$. The theta function is
holomorphic and quasi-periodic with respect to the lattice $\Bbb
Z^g + \Omega \Bbb Z^g \subset \Bbb C^g$ and it is either even or
odd according to the parity of $\xi$.  By definition the parity
of $\xi$ is the parity of $q(\xi) ={}^t a\cdot b$, and
$\vartheta[\xi](-z, \Omega) = (-1)^{q(\xi)} \vartheta[\xi](z,
\Omega)$.  In particular, if $\Theta[\xi] \subset A = \Bbb
C^g/(\Bbb Z^g + \Omega \Bbb Z^g)$ denotes the zero divisor of
$\vartheta[\xi](z, \Omega)$, then $\Theta[\xi]$ is a symmetric
theta divisor, and $\text{mult}_0(\Theta[\xi])\equiv q(\xi) \pmod
2$.

The more difficult part of the theory describes the
transformation  of $\vartheta[\xi](z, \Omega)$ under the natural
action of $\text{Sp}(2g, \Bbb Z)$ on $\{0, 1\}^{2g} \times \Bbb
C^g \times \Cal H_g$.  Our notation for the action of $\gamma \in
\text{Sp}(2g, \Bbb Z)$ on $(\xi, z, \Omega)$ is $(\hat{\xi},
\hat{z}, \hat{\Omega}) = \gamma \cdot (\xi, z, \Omega)$; the
formula for $[\hat{\xi}] = \gamma \cdot [\xi]$ is given in
\cite{I1, p\. 226}.  The transformation law then relates
$\vartheta[\hat{\xi}](\hat{z}, \hat{\Omega})$ to
$\vartheta[\xi](z, \Omega)$; i.e., there exists a function
$\sigma = \sigma(\gamma, \xi, z, \Omega)$ such that $$ \align
\vartheta[\hat{\xi}](\hat{z}, \hat{\Omega}) &= \sigma(\gamma,
\xi, z, \Omega) \vartheta[\xi](z, \Omega),\\ \sigma(\gamma, \xi,
z, \Omega) &= \kappa(\gamma)e^{2\pi i
\phi_{\xi}(\gamma)}\det(C\Omega + D)^{1/2}e^{\pi i Q_{\gamma,
\Omega}(z)}. \endalign $$ (There may be a choice of sign in
$\det(C\Omega + D)^{1/2}$ and in $\kappa(\gamma)$ separately, but
the product $\kappa(\gamma)\det(C\Omega + D)^{1/2}$ is a
well-defined function of $(\gamma, \Omega)$.  For the properties
of $\kappa(\gamma)$ and the formula for $\phi_{\xi}(\gamma)$, see
\cite{I1}, \cite{I2}; the formula for $Q_{\gamma, \Omega}(z)$
will be given below.)

For any theta characteristic $\xi$, define $$ \Gamma[\xi] =
\{\gamma \in \text{Sp}(2g, \Bbb Z)\ |\ \gamma \cdot [\xi] = [\xi]
\}. $$ Then in particular, putting $\sigma_{\xi}(\gamma, z,
\Omega) = \sigma(\gamma, \xi, z, \Omega)$, there exists an
automorphy factor $\sigma_{\xi}: \Gamma[\xi] \times (\Bbb C^g
\times \Cal H_g) \rightarrow \Bbb C^*$ for the action of
$\Gamma[\xi]$ on $\Bbb C^g \times \Cal H_g$ such that
$\vartheta[\xi](\hat{z}, \hat{\Omega}) = \sigma_{\xi}(\gamma, z,
\Omega) \vartheta[\xi](z, \Omega)$ for every $\gamma \in
\Gamma[\xi]$ and $(z, \Omega) \in \Bbb C^g \times \Cal H_g$.

Note that $\Gamma(2) \subset \Gamma[\xi]$ for each theta
characteristic $\xi$.  In fact, if $q_{\xi}$ denotes the $\Bbb
Z/2$-quadratic form corresponding to the theta characteristic
$\xi$ and $O(q_{\xi}) \subset \text{Sp}(2g, \Bbb Z/2)$ denotes
the orthogonal group defined by $q_{\xi}$, then $\Gamma[\xi]$ is
the preimage of $O(q_{\xi})$ under the canonical surjection
$\text{Sp}(2g, \Bbb Z) \rightarrow \text{Sp}(2g, \Bbb Z/2)$.
Since $\text{Sp}(2g, \Bbb Z)$ acts transitively on the set of odd
theta characteristics \cite{I3, Cor\. p\. 213}, the subgroups
$\Gamma[\xi]$, $\xi$ odd, are all conjugate in $\text{Sp}(2g,
\Bbb Z)$, and for any one of them, $\Cal H_g/\Gamma[\xi] \cong
\tilde{\Cal A}_g$.

To render the transformation law for $\vartheta[\xi]$ as simple
as possible, we pass to the important subgroup $$\multline
\Gamma' = \Gamma(4, 8)=\Big\{\left(\smallmatrix A & B \\ C & D
\endsmallmatrix\right)\in\text{Sp}(2g,\Bbb Z)\ \ |\
\left(\smallmatrix A & B \\ C & D \endsmallmatrix\right) \equiv
I_{2g} \pmod4, \text{ and } \\ \text{diagonal entries of } A\,^tB
\text{ and } C\,^tD \equiv 0\pmod8 \Big\} \endmultline$$  (from
\cite{I1, p\. 220}); it lies between the congruence subgroups
$\Gamma(8)$ and $\Gamma(4)$ and is normal in $\text{Sp}(2g,\Bbb
Z)$ \cite{I1, p\. 222}.  Thus we have the containments $$
\Gamma(8) \subset \Gamma(4, 8) \subset \Gamma(4) \subset
\Gamma(2) \subset \Gamma[\xi] \subset \text{Sp}(2g, \Bbb Z). $$

\noindent 1.4.1 \ {\sl The transformation law for
$\vartheta[\xi](z,\Omega)$ with respect to $\Gamma'$.} For
$\gamma \in \Gamma'$ and $(z, \Omega) \in \Bbb C^g \times \Cal
H_g$, if we write $$\spreadlines{-.04in}\align (\hat{z},
\hat{\Omega}) = \gamma\cdot (z, \Omega) &= ({}^t(C\Omega +
D)^{-1}z, (A\Omega + B)(C\Omega + D)^{-1}),\\ \intertext{then}
\vartheta[\xi](\hat{z}, \hat{\Omega}) &= \det (C\Omega +D)^{1/2}
e^{\pi i Q_{\gamma, \Omega}(z)} \vartheta[\xi](z, \Omega),
\endalign $$ where $ Q_{\gamma, \Omega}(z) = \ ^t z(C\Omega +
D)^{-1} Cz$ \cite{I1, p\. 227}.

That is, for $\gamma \in \Gamma'$ there is a canonical choice of
square root of $\det(C\Omega + D)$ so that $\sigma(\gamma, z,
\Omega) = \det(C\Omega + D)^{1/2} e^{\pi i  Q_{\gamma,
\Omega}(z)}$ is the transformation factor, which is now
independent of $\xi$.  We let $\rho: \Gamma' \times \Cal H_g
\rightarrow \Bbb C^*$ denote the fundamental automorphy factor
$\rho(\gamma, \Omega) = \det(C\Omega + D)^{1/2}$.

For us the main point will be following.  For any fixed odd
characteristic $\xi$, consider the Taylor expansion of the theta
function $\vartheta[\xi](z, \Omega)$ with respect to $z$ about
the origin: $$ \vartheta[\xi](z, \Omega) = \ell(z, \Omega) + m(z,
\Omega) + \cdots , $$ where $\ell = \ell_{\xi}$, $m = m_{\xi}$,
$\dots$, are respectively the  linear, cubic and higher odd order
terms in the variable $z$.  Then the transformation law (1.4.1)
shows that all of the ideals $(\ell)$, $(\ell, m)$, $\dots$,
transform by the same automorphy factor $\rho$. In other words,
by general principles on jet prolongation of sections of bundles,
the automorphy factor $\rho$ governs every homogeneous Taylor
term of $\vartheta[\xi]$ modulo all the previous terms.

Let $P$ denote the fundamental (positive) line bundle on $\Cal
A'_g$ with automorphy factor $\rho$.  As we see from the
transformation law, the classical theta nulls $\vartheta[\xi](0,
\Omega)$ for {\it even\/} $\xi$, define sections of $P$ (and an
embedding theorem is proved in \cite{I1}, \cite{I2}, \cite{I3};
cf\. \cite {M2, p\. 73}), but we will employ this line bundle in
connection with the {\it odd\/} theta functions. Clearly, $P^2
\cong \det H = \Lambda$, since $\rho^2 = \det(C\Omega + D)$ and
$\eta(\gamma, \Omega) = C\Omega + D$ is the automorphy factor for
$H$.  Let $N = P^*$ be the fundamental (negative) line bundle on
$\Cal A'_g$ with automorphy factor $\nu(\gamma, \Omega) =
\det(C\Omega + D)^{-1/2}$ with respect to $\Gamma'$.

\remark{ \rm{1.4.2}\ Remark} Apparently, the line bundle $N$
corresponds to Kempf's $Q(A)$ in \cite{K2, p\. 70}, and his Cor\.
8.7 \cite{K2, p\. 72} reflects the progressively simpler
transformation laws for $\vartheta$. \endremark

\vskip.5\baselineskip \subhead 1.5\quad The homomorphism
$\ell_{\xi}$ \endsubhead \vskip.5\baselineskip

Now we consider the linear term (in $z$) of an odd theta function
$\vartheta[\xi](z, \Omega)$ at the origin of $\Bbb C^g$, as
$\Omega$ varies over $\Cal H_g$.  We will show that this linear
term can be viewed as a homomorphism $N \rightarrow H$ over $\Cal
A'_g$, with  the geometric interpretation that for each $y = (A,
\Theta, \xi') \in \Cal A'_g$, the image of the fibre $N_y$ in
$H_y = T_0^*A$ is the conormal line at $0$ to $\Theta[\xi]$ in
$A$ when the theta divisor $\Theta[\xi]$ is nonsingular at $0$
(and the image is $0$ when $\Theta[\xi]$ is singular at $0$).

\proclaim{1.5.1 Proposition}  For each odd theta characteristic
$\xi$, $\ell(z, \Omega) = \nabla_z|_0\big(\vartheta[\xi](z,
\Omega)\big)$ defines a section of $P \otimes H$ over $\Cal
A'_g$, or equivalently, a homomorphism $$ \ell_{\xi}: N
\rightarrow H $$ of bundles on $\Cal A'_g$. \endproclaim

\demo{Proof}  Here $\nabla_z|_0\big(\vartheta[\xi](z,
\Omega)\big)$ is the  column vector, with entries $a_j(\Omega) =
\d{z_j}|_0\big(\vartheta[\xi](z, \Omega)\big)$, $1 \le j \le g$,
which represents the linear form $\sum a_j(\Omega)z_j$, the
differential of $\vartheta[\xi](z, \Omega)$ at $z = 0$.  We check
that $\Omega \mapsto \nabla_z|_0\big(\vartheta[\xi](z,
\Omega)\big)$, $\Omega \in \Cal H_g$, represents a section of $P
\otimes H$.  First, from the transformation law (1.4.1) we deduce
that $$ \vartheta[\xi](z, \gamma\cdot\Omega) = \det(C\Omega +
D)^{1/2} e^{\pi i \tilde{ Q}_{\gamma,
\Omega}(z)}\vartheta[\xi](^t(C\Omega + D)z, \Omega), $$ where
$\tilde{Q}_{\gamma, \Omega}(z) = Q_{\gamma, \Omega}(^t(C\Omega +
D)z) ={}^t z\,C\,^t(C\Omega + D)z$ is a homogeneous quadratic
function of $z$.

Now we compute $\nabla_z|_0\big(\vartheta[\xi](z, \gamma
\cdot\Omega)\big)$:  $$\align \d{z_j}\big(\vartheta[\xi](z,
\gamma \cdot \Omega)\big) &= \d{z_j}\left(\det(C\Omega +
D)^{1/2}e^{\pi i\tilde{Q}_{\gamma,
\Omega}(z)}\vartheta[\xi](^t(C\Omega + D)z, \Omega)\right)\\
&=\det(C\Omega + D)^{1/2}\d{z_j}\left((1 + h_{\gamma,
\Omega}(z))\vartheta[\xi](^t(C\Omega + D)z, \Omega)\right)\\
&=\det(C\Omega + D)^{1/2}\left(\d{z_j}(h_{\gamma,
\Omega}(z))\vartheta[\xi](^t(C\Omega + D)z, \Omega)) + \right.\\
&\hskip.5in\left. (1 + h_{\gamma,\Omega}(z))
\d{z_j}\big(\vartheta[\xi](^t(C\Omega + D)z, \Omega)\big)\right),
\endalign $$  where $h_{\gamma, \Omega}$ has order $2$ and higher
in $z$.  Thus $$ \nabla_z|_0\big(\vartheta[\xi](z, \gamma \cdot
\Omega)\big) = \det(C\Omega +
D)^{1/2}\nabla_z|_0\left(\vartheta[\xi]\big(^t(C\Omega + D)z,
\Omega\big)\right). $$ By the chain rule, $$ \align
\nabla_z|_0\left(\vartheta[\xi]\big(^t(C\Omega + D)z,
\Omega\big)\right) &= \eta^{(1)}(\gamma,
\Omega)\left(\nabla_z|_0\big(\vartheta[\xi](z,
\Omega)\right)\big)\\ &= (C\Omega +
D)\nabla_z|_0\big(\vartheta[\xi](z, \Omega)\big), \endalign $$
where $\eta^{(1)}(\gamma, \Omega)$ is the appropriate cotangent
map.  Therefore, $$ \nabla_z|_0\big(\vartheta[\xi](z, \gamma
\cdot \Omega)\big) = \det (C\Omega +
D)^{1/2}(C\Omega+D)\nabla_z|_0
\big(\vartheta[\xi](z,\Omega)\big). $$  Thus we have shown that
$\ell(z, \Omega) = \nabla_z|_0\big(\vartheta [\xi](z,
\Omega)\big)$ satisfies the transformation rule $$ \ell(z, \gamma
\cdot \Omega) = \det(C\Omega + D)^{1/2} (C\Omega + D) \ell(z,
\Omega), $$ for all $\gamma \in \Gamma'$ and $\Omega \in \Cal
H_g$, where $\det(C\Omega + D)^{1/2}(C\Omega + D)$ is the
automorphy factor for $P \otimes H$.  Q.E.D. \enddemo

For $\gamma \in \text{Sp}(2g, \Bbb Z)$, we may regard the map $L
= L(\gamma, \Omega): \Bbb C^g \rightarrow \Bbb C^g, z \mapsto \
^t(C\Omega + D)^{-1}z$, as an isomorphism from the universal
cover of the abelian variety $A_{\Omega} = \Bbb C^g/(\Bbb Z^g +
\Omega \Bbb Z^g)$ to that of the abelian variety
$A_{\hat{\Omega}} = \Bbb C^g/(\Bbb Z^g + \hat{\Omega}\Bbb Z^g)$,
where $\hat{\Omega} = \gamma \cdot \Omega$. Then the derivative
$\eta_{(1)}(\gamma, \Omega)$ of $L(\gamma, \Omega)$ at $0$ is the
map $z \mapsto \mu(\gamma, \Omega)z$ given by the automorphy
factor $\mu(\gamma, \Omega) ={}^t(C\Omega + D)^{-1}$, and the
cotangent map $\eta^{(1)}(\gamma, \Omega)$ of $L(\gamma, \Omega)$
at $0$ is the dual map $z \mapsto \eta(\gamma, \Omega)z$ given by
the automorphy factor $\eta(\gamma, \Omega) = C\Omega + D$.

Note that, although the homomorphism $\ell_{\xi}: N \rightarrow
H$ is injective as a map on sheaves of sections, it induces the
$0$-map on the fibres on $N$ over points of $\Cal A'_g$ for which
the theta divisor $\Theta[\xi]$ is singular at the origin.

\head 2. The cubic hypersurface associated to an odd theta
characteristic \endhead

\subhead 2.1\quad The open set $\Cal U_{\xi}$ \endsubhead
\vskip.5\baselineskip

Recall that there exists, in addition to a universal family
$\frak A \rightarrow \Cal A'_g$ of abelian varieties, a universal
family $\sctheta \rightarrow \Cal A'_g$ of theta divisors
$\{\Theta_y = \Theta[\xi] \subset A_y\ |\ y = (A, \Theta, \xi')
\in \Cal A'_g\}$.  In particular, for each $(A, \Theta, \xi') \in
\Cal A'_g$, the marked symmetric theta divisor $\Theta[\xi]
\subset A$ has odd multiplicity at the origin $0$ of $A$.  Since
we want to look at the tangent hyperplanes of these theta
divisors at the origin, we will pass to the open subset of $\Cal
A'_g$ over which the theta divisors have multiplicity exactly one
at the origin.  Let $$ \Cal U_{\xi} = \{y \in \Cal A'_g\ |\
\Theta_y \text{ is nonsingular at } 0\}. $$ Then $\Cal U_{\xi}
\subset \Cal A'_g$ is an open subset and on $\Cal U = \Cal
U_{\xi}$ we have a rank $g - 1$ vector bundle $\Cal E = \Cal
E_{\xi}$ whose fibre at each point $y\in \Cal U$ is $E_y =
T_0\Theta_y$.  Note that $\Cal E$ is naturally a subbundle of the
restriction to $\Cal U$ of the rank $g$ vector bundle $\Cal V =
H^*$ whose fibre at $y$  is $T_0(A)$.  Recall that on $\Cal U$ we
have the line subbundle $\Cal L_{\xi} = \ell_{\xi}(N)$ of $H =
\Cal V^*$ defined by the conormal lines to the theta divisors
$\Theta[\xi]$ at the origin of the abelian varieties.  Thus, by
construction, $\Cal E$ is the rank $g - 1$ subbundle $\Cal
L_{\xi}^{\perp}$ of $\Cal V$ over the open set $\Cal U$ and we
have the exact sequence $$ 0 \rightarrow \Cal L_{\xi} \rightarrow
H \rightarrow \Cal E^*_{\xi} \rightarrow 0 $$ of vector bundles
on $\Cal U$.

The next step is to use the cubic term (in $z$) of the theta
function $\vartheta[\xi](z, \Omega)$ to construct cubic forms on
the fibres of the vector bundle $\Cal E$ over $\Cal U$.  It will
be quite important in our applications that in passing from $\Cal
A'_g$ to $\Cal U$ we have not lost information about divisors of
$\Cal A'_g$.

\proclaim{2.1.1 Lemma} For each $g \ge 4$, $Z = \Cal A'_g - \Cal
U$ is a closed algebraic subset of $\Cal A'_g$ of codimension
$\ge 2$. \endproclaim

\demo{Proof} First we check that $Z \subset \Cal A'_g$ is
algebraic. Now $\sctheta \rightarrow \Cal A'_g$ is a proper
morphism of algebraic varieties. The critical locus $\Sigma
\subset \sctheta$ is an algebraic subset of $\sctheta$ with the
property that, for each $y = (A, \Theta, \xi') \in \Cal A'_g$,
the fibre of $\Sigma$ over $y$ is $\Sigma_y =
\text{Sing}(\Theta_y)$, the singular locus of the theta divisor
$\Theta_y \subset A_y$.  Also the zero section $E$ of the family
of abelian varieties $\frak A \rightarrow \Cal A'_g$ is an
algebraic subset of $\sctheta$.  Therefore, the intersection $E
\cap \Sigma = \{(y, 0)\ |\ 0 \text{ is singular on }\Theta_y\}
\subset \sctheta$ is an algebraic subset of $E$, and hence under
the isomorphism $E \rightarrow \Cal A'_g$, the image $Z = \{y\ |\
0 \text{ is singular on } \Theta_y\}$ of $E \cap \Sigma$ is
algebraic in $\Cal A'_g$. \enddemo

Now we check that the codimension of $Z$ in $\Cal A'_g$ is at
least 2. Denote by $\script J'_g$ the preimage in $\Cal A'_g$ of
the Jacobian locus $\script J_g\subset\Cal A_g$. Using the fact
that for any $g \ge 4$ the Picard number of $\Cal A'_g$ is equal
to one, it suffices to check that $Z \cap \script J'_g$ has
codimension $\ge 2$ in $\script J'_g$. (Cf\. \cite{S-V, Cor\.
(0.7), p\. 350}.) (For $g \ge 4$ and for any arithmetic subgroup
$\Gamma$ of $\text{Sp}(2g, \Bbb R)$, $\text{Pic}(\Cal H_g/\Gamma)
\otimes \Bbb Q \cong \Bbb Q$,  by the result of \cite{Bo} that
$H^2(\Gamma, \Bbb Q) \cong \Bbb Q$ (cf\. \cite{M4, pp\. 354-355},
\cite{S-V, Part 0, \S C}).  For our purposes, it suffices to work
in $\Cal A_g$ and use the result of \cite{Fr1} that for $g \ge
3$, $H^2(\text{Sp}(2g, \Bbb Z), \Bbb Z) \cong \Bbb Z$.)  Now we
pass to level $1$ moduli, i.e., from $Z \cap \script J'_g\subset
\Cal A'_g$ to its image $\script J_g \subset \Cal A_g$. Consider
the locus of Jacobians which have at least one odd theta
characteristic $\xi$ such that the corresponding theta divisor
$\Theta[\xi]$ is singular at $0$; this type of locus has been
studied by Harris and Teixidor. In the moduli space of curves
$\Cal M_g$, the locus in question is denoted $\Cal M_g^2$; it is
the locus of genus $g$ curves having an odd theta characteristic
with $h^0 \ge 3$.  From \cite{T, Thm. (2.13), p\. 109} we see
that every component of $\Cal M_g^2$ has codimension at least $3$
in $\Cal M_g$. Q.E.D.

\remark{\rm {2.1.2}\  Remark}  Here are two possible alternative
arguments, both based on the use of the codimension one boundary
of the Igusa-Mumford compactification $\bar{\Cal A}_g$ of $\Cal
A_g$ (as in \cite{M4}, \cite{D, Lemme 2.1, p\. 698}).

(1) R.~Salvati~Manni (oral communication) has calculated that the
closure of $Z$ in $\bar{\Cal A}_g$ has the expected codimension
$g$ on the boundary; in particular, there is no codimension 1
component of $Z$.

(2) By definition, $Z \subset \Cal N_0 = \{(A, \Theta, \xi') \in
\Cal A'_g\ |\ \Theta \text{ is singular}\}$. Debarre \cite{D} has
shown that in the level $1$ moduli space, $\Cal N_0$ has exactly
two irreducible components, both of codimension one, and over the
generic point of each of these two irreducible components, the
only singularities of the theta divisor are double points.  On
the other hand, over every point $y \in Z$, there is a
singularity of $\Theta_y$ that is not a double point, and hence
$Z$ cannot contain an irreducible component of $\Cal N_0$.
Therefore, $\text{codim}(Z) > \text{codim}(\Cal N_0) = 1$.
\endremark

\proclaim{2.1.3 Corollary}  For any holomorphic line bundle $\Cal
M$ on $\Cal A'_g$, $H^0(\Cal A'_g, \Cal M) \rightarrow H^0(\Cal
U, \Cal M)$ is an isomorphism. \endproclaim

\demo{Proof}  It suffices to establish the isomorphism of
sections locally on $\Cal A'_g$; since the line bundle $\Cal M$
is locally trivial on $\Cal A'_g$, it is enough to know that for
all sufficiently small open sets $U \subset \Cal A'_g$, $H^0(U,
\Cal O) \rightarrow H^0(U - Z, \Cal O)$ is an isomorphism.  Since
$Z$ has complex codimension $\ge 2$, $\Cal O(U) \rightarrow \Cal
O(U - Z)$ is an isomorphism by Hartogs' Theorem \cite{G-R, p\.
132}.  Q.E.D. \enddemo

\remark{\rm{2.1.4} \ Remark}  Rather than simply removing $Z$ to
define the bundle $\Cal E$, one can utilize a natural Nash
blow-up construction that yields a (proper, birational)
modification $\sigma: \Cal Y \rightarrow \Cal A'_g$ such that the
bundle inclusion $\Cal L = \Cal L_{\xi} \subset H$ on $\Cal U$
extends to an inclusion $\hat{\Cal L} \subset \sigma^*(H)$ on
$\Cal Y$, and hence the vector bundle $\Cal E$ on $\Cal U \subset
\Cal A'_g$ extends to a vector subbundle $\hat{\Cal E} =
\hat{\Cal L}^{\perp} \subset \sigma^*(\Cal V)$ on all of $\Cal
Y$. \endremark

\vskip.5\baselineskip \subhead 2.2\quad The homomorphism
$\overline{m}_{\xi}$ \endsubhead \vskip.5\baselineskip

We now consider the cubic Taylor term of an odd theta function
$\vartheta[\xi](z, \Omega)$, modulo the linear term.

\proclaim{2.2.1 Proposition}  For each odd theta characteristic
$\xi$, $$
\tfrac1{3!}\nabla_z^3|_0\big(\vartheta[\xi](z,\Omega)\big)
\quad\left(\text{\rm mod }\nabla_z|_0\big(\vartheta[\xi](z,
\Omega)\big)\right) $$ defines a section of $P\otimes S^3\Cal
E_{\xi}^*$ over $\Cal U \subset \Cal A'_g$, or, equivalently, a
homomorphism $$ m_{\xi}: N \rightarrow S^3\Cal E_{\xi}^* $$ of
bundles on $\Cal U$. \endproclaim

\demo{Proof}  Here $\nabla_z^3|_0\big(\vartheta[\xi](z,
\Omega)\big)$ is the cubic form $\sum_{1\le j,k,l\le g}
a_{jkl}(\Omega)z_j z_k z_l$, where $a_{jkl}(\Omega) =$\break
$\frac{\partial^3}{\partial z_j\partial z_k\partial
z_l}|_0\big(\vartheta[\xi](z, \Omega)\big)$.  Let $p: \Cal H_g
\rightarrow \Cal H_g/\Gamma' = \Cal A'_g$ be the quotient map,
and assume $\Omega\in p^{-1}(\Cal U)$. We will check that
$$\Omega\mapsto \nabla_z^3|_0\big(\vartheta[\xi](z, \Omega)\big)
\quad\left(\text{mod }\nabla_z|_0\big(\vartheta[\xi](z,
\Omega)\big)\right), $$ represents a section of $P \otimes
S^3\Cal E_{\xi}^*$ over $\Cal U$. Consider the Taylor expansion
$\vartheta[\xi](z, \Omega) = \ell(z, \Omega) + m(z, \Omega) +
\dots$. We apply the transformation law (1.4.1), expanding
$e^{\pi i Q_{\gamma,\Omega}(z)} = 1 + q_{\gamma, \Omega}(z) +
\dots$, where $q_{\gamma, \Omega}(z) = \pi i  Q_{\gamma,
\Omega}(z)$ is homogeneous quadratic in $z$.  Then $$ \align
\ell(\hat{z}, \hat{\Omega}) + m(\hat{z}, \hat{\Omega}) + \dots &=
\det(C\Omega + D)^{1/2}(1 + q_{\gamma, \Omega}(z) + \dots
0)(\ell(z, \Omega) + m(z, \Omega) + \dots )\\ &= \det(C\Omega +
D)^{1/2}\left(\ell(z, \Omega) + \big(m(z, \Omega) + q_{\gamma,
\Omega}(z) \ell(z, \Omega)\big) + \dots\right), \endalign $$ so
we obtain $\ell(\hat{z}, \hat{\Omega}) = \det(C\Omega +
D)^{1/2}\ell(z, \Omega)$ (as before) and $$ m(\hat{z},
\hat{\Omega}) = \det(C\Omega + D)^{1/2}\big(m(z, \Omega) +
q_{\gamma, \Omega}(z)\ell(z, \Omega)\big). $$ Now we replace
$\hat{z}$ on the left by $z$ and hence $z$ on the right by
$\tilde{z} ={}^t(C\Omega + D)z$.  Then we have $$ m(z, \gamma
\cdot \Omega) = \det(C\Omega + D)^{1/2} (m(\tilde{z}, \Omega) +
q_{\gamma, \Omega}(\tilde{z})\ell(\tilde{z}, \Omega)). $$ In
other words, $$ m(z, \gamma \cdot \Omega) = \det(C\Omega +
D)^{1/2} \eta^{(3)}(\gamma, \Omega)\big(m(z, \Omega) + q_{\gamma,
\Omega}(z)\ell(z, \Omega))\big), $$ where $\eta^{(3)}(\gamma,
\Omega)$ is the natural map on cubic forms on $\Bbb C^g$ induced
by the cotangent map $\eta^{(1)}(\gamma, \Omega)$.

Now we consider the cubic form $m(z,\Omega)$ modulo multiples of
$\ell(z, \Omega)$ (by quadratic forms). Equivalently, we restrict
the cubic form to the hyperplane $\{z \in \Bbb C^g\ |\ \ell(z,
\Omega) = 0\}$.  Note that this hyperplane is exactly the fibre
$E_{\Omega} \subset \Cal V_{\Omega} = \Bbb C^g$ of $p^*(\Cal E)$
at $\Omega$.  Thus, the resulting cubic form $\overline{m}(z,
\Omega)$ on the fibres of $p^*(\Cal E)$ transforms as follows: $$
\overline{m}(z, \gamma \cdot \Omega) = \det(C\Omega + D)^{1/2}
\overline{\eta}^{(3)}(\gamma, \Omega)(\overline{m}(z, \Omega)),
$$ where $\overline{\eta}^{(3)}(\gamma, \Omega)$ is the induced
map on cubic forms $\text{mod } \ell(z, \Omega)$.  This proves
that $\overline{m}(z, \Omega)$ transforms as a section of $P
\otimes S^3\Cal E_{\xi}^*$.  Q.E.D. \enddemo

\remark{\rm{2.2.2}\  Remark}  After the modification $\sigma:
\Cal Y \rightarrow \Cal A'_g$ indicated in Remark (2.1.4), the
homomorphism $\overline{m}: N \rightarrow S^3\Cal E^*$ over $\Cal
U \subset \Cal A'_g$ extends to a homomorphism $\widehat{m}:
\sigma^*(N) \rightarrow S^3\hat{\Cal E}^*$ over all of $\Cal Y$.
\endremark \vskip.5\baselineskip

Thus, for each point $y = (A, \Theta, \xi') \in \Cal U \subset
\Cal A'_g$, if the cubic form $\overline{m}_{\xi}$ is nonzero,
then it defines a distinguished cubic hypersurface $V$ in $\Bbb
E_{\xi} = \Bbb P T_0(\Theta[\xi]) \cong \Bbb P^{g-2}$, and $V
\subset \Bbb P^{g-2}$ depends (up to projective equivalence) only
on the image point $(A, \Theta, \xi) \in \tilde{\Cal A}_g$.  For
$(A, \Theta, \xi) \in \tilde{\Cal A}_g$, let $V_\xi$ denote the
cubic hypersurface of $\Bbb E_{\xi} = \Bbb P(\ell_{\xi} = 0)$
defined in $\Bbb P^{g-1}$ by the ideal $(\ell_{\xi}, m_{\xi})$
(provided that $\ell_{\xi}, m_{\xi}$ is a regular sequence).

\vskip.5\baselineskip \subhead 2.3\quad Jacobians \endsubhead
\vskip.5\baselineskip

Recall the correspondence between the theta characteristics of a
curve and those of its Jacobian.  A {\it theta characteristic\/}
of a genus $g$ curve $C$ is a degree $g - 1$  line bundle $\xi$
(up to isomorphism) such that $\xi^{\otimes 2} \cong \Omega_C$,
the canonical bundle of $C$.  The {\it parity\/} of the theta
characteristic $\xi$ is the parity of $h^0(\xi)$. The one-to-one
correspondence (preserving parity) between the theta
characteristics of a curve $C$ and those of its Jacobian can be
found in \cite{ACGH, pp\. 281-294}, \cite{C, pp\. 140-143},
\cite{R-F, pp\. 176-177}, \cite{M3, pp\. 162-170}. Let $\Cal M_g$
be the moduli space of genus $g$ curves over $\Bbb C$ and let
$\script J_g$ be its image in $\Cal A_g$.  Thus $\Cal M_g$ is the
set of isomorphism classes of smooth, connected, complete curves
of genus $g$ over $\Bbb C$, and $\script J_g$ is the
(Zariski-locally-closed) subset of $\Cal A_g$ consisting of
isomorphism classes of principally polarized Jacobians of genus
$g$ curves.

An odd theta characteristic $\xi$ of a genus $g$ curve $C$ is
{\it nondegenerate\/} if the complete linear series $|\xi| =
\{D\}$ for a single degree $g - 1$ divisor $D = p_1 + \dots +
p_{g-1}$ consisting of $g - 1$ distinct points $p_1, \dots,
p_{g-1}$.  Recall that an odd theta characteristic $\xi$ is {\it
nonsingular\/} if $h^0(\xi) = 1$; i.e., $|\xi| = \{D\}$ for a
single effective degree $g - 1$ divisor $D$.  (By the Riemann
singularities theorem, $h^0(\xi) = \text{mult}_0(\Theta[\xi])$,
so an odd theta characterisic $\xi$ of $C$ is nonsingular if and
only if the origin $0$ is a nonsingular point of $\Theta[\xi]$.)
Note that by geometric Riemann-Roch (cf\. \cite{ACGH, p\. 12}),
if $D$ is any effective divisor of degree $g - 1$ then $h^0(\Cal
O(D)) = 1$ if and only if $D$ spans a hyperplane of $\Bbb
P^{g-1}$.

\example{2.3.1 Examples}  (i)  Let $C \subset \Bbb P^2$ be a
nonsingular plane quartic $(g = 3)$ with a higher flex $p$ (i.e.,
a weight $2$ Weierstrass point).  Then $\xi = \Cal O(2p)$ is a
degenerate nonsingular odd theta characteristic.  Indeed, $|\xi|
= |K - 2p|$ is cut out by the lines that are tangent to $C$ at
$p$, and there is only one such line.

\noindent (ii)  Let $C$ be a hyperelliptic curve of genus $g = 4k
+ 1$, $k \ge 1$, and let $\Gamma \subset \Bbb P^{4k}$ be its
canonical image. Then $\xi = \Cal O(2k \cdot g^1_2)$ is a
singular odd theta characteristic. For if $D$ consists of the sum
of the $g^1_2$ divisors over $2k$ of the branch points on
$\Gamma$, then $|K - \xi| = |K - D|$ is cut out by the
hyperplanes through the $2k$ points in $\Bbb P^{4k}$, so this
linear system has projective dimension $2k$; hence $h^0(\xi) =
h^0(K - \xi) = 2k + 1$.

\noindent (iii)  In genus $4$, all the odd theta characteristics
of a hyperelliptic curve are nondegenerate.  By choosing any $3$
of the $10$ Weierstrass points, we get $120$ nondegenerate odd
theta characteristics, so these are all the odd theta
characteristics. \endexample

\proclaim{2.3.2 Proposition}  For $g \ge 3$, let $N\Cal M_g = \{C
\in \Cal M_g\ |\ $ every odd theta characteristic of $C$ is
nondegenerate\} and let $N \script J_g$ be the image of $N\Cal
M_g$ in $\script J_g$. \roster \item"{(i)}" $N\Cal M_g$ is a
nonempty Zariski-open subset of $\Cal M_g$, and hence $N\script
J_g$ is Zariski-dense in $\script J_g$. \item"{(ii)}"  For each
$(J, \Theta) \in N\script J_g$ and each odd theta characteristic
$\xi$ of $(J, \Theta)$, let $\Omega \in \Cal H_g$ be a period
matrix for $(J, \Theta)$, and let $\ell_{\xi}$ and $m_{\xi}$ be
the linear and cubic terms, respectively, of the theta function
$\vartheta[\xi](z, \Omega)$.  Let $\Bbb E_{\xi}$ be the
hyperplane of canonical space $\Bbb P^{g-1} = \Bbb P T_0(J)$
defined by the linear form $\ell_{\xi}$, and let $V_\xi$ be the
cubic hypersurface of $\Bbb E_{\xi}$ defined in $\Bbb P^{g-1}$ by
the ideal $(\ell_{\xi}, m_{\xi})$. Then $V_\xi$ is a Fermat
cubic. \endroster \endproclaim

\demo{Proof of (i)}  First we prove that $N\Cal M_g$ is a
Zariski-open subset of $\Cal       M_g$. Let $\{C\}$ be an
irreducible family of genus $g$ curves that dominates $\Cal M_g$
by a finite proper map. Let $S=\{C\ |\ \text{some odd theta
characteristic of }C\text{ is singular}\}$, and let $D=\{C\ |\
\text{some odd theta}$ $\text{characteristic of }C\text{ is
degenerate}\}$. By upper semicontinuity of $h^0$, $S$ is
Zariski-closed, so it suffices to show that $D-S$ is
Zariski-closed in the complement of $S$. To show this we pass to
the family $\{(C,\xi)\}$, where $\xi$ is a nonsingular odd theta
characteristic, so that $|\xi|$ has a unique representative
$p_1+\cdots + p_{g-1} \in C^{(g-1)}$, and use the fact that the
discriminant is Zariski-closed in the total space of the family
$\{C^{(g-1)}\}$.

Next we show that $N\Cal M_g$ is nonempty. (1) There exists a
genus $g$ curve with a nondegenerate odd theta characteristic.
For example, a genus $g$ hyperelliptic curve $C$ has $2g+2$
distinct Weierstrass points $p_1, \dots, p_{2g+2}$, and if $D =
p_1+\dots+p_{g-1}$, then $\xi=\Cal O(D)$ is a nondegenerate odd
theta characteristic. (2) There exists a family of genus $g$
curves over a connected, smooth base $B$ such that $B$ dominates
$\Cal M_g$ and the monodromy of the family is transitive on the
set of odd theta characteristics. An example is the universal
family of curves over the nonempty open subset of $\Cal M_g$ of
automorphism-free curves (cf\. \cite{ACGH, p\. 294}). Now (1) and
(2) imply that all the odd theta characteristics of a generic
genus $g$ curve are nondegenerate.

Finally, since $N\Cal M_g$ is nonempty Zariski-open in $\Cal M_g$
and $\Cal M_g$ is irreducible, $N\Cal M_g$ is Zariski-dense in
$\Cal M_g$, so its image $N\script J_g$ in $\script J_g$ is
Zariski-dense under the surjection $\Cal M_g\to \script J_g$.
Q.E.D.\enddemo

\remark{\rm{2.3.3}\  Remark}
That $N\Cal M_g$ is nonempty is a special case of a general
result of J.~McKernan \cite{Mc} on the structure of the loci of
hyperplane sections of a generic canonical curve that have a
given type (set of multiplicities); in our case the type is $(2,
\dots, 2)$. \endremark

\demo{Proof of (ii)}  For $(J, \Theta) \in N\script J_g$ and any
odd theta characteristic $\xi$ of $(J, \Theta)$, let $C$ be a
genus $g$ curve with (polarized) Jacobian $J(C) \cong (J,
\Theta)$, such that $\xi$ is a nondegenerate odd theta
characteristic of $C$.  We will use $C$ to determine the cubic
hypersurface $V_\xi$.  Let $D = p_1 + \dots + p_{g-1}$ be the
unique divisor of $|\xi|$ on $C$.  The points $p_1$, $\dots$,
$p_{g-1}$ are distinct and their images $\varphi(p_1)$, $\dots$,
$\varphi(p_{g-1})$ under the canonical map $\varphi: C
\rightarrow \Bbb P^{g-1}$ span a (unique) hyperplane $H$. In
particular, the points $\varphi(p_i)$, $i = 1,\dots,g-1$, are
linearly independent, and at each of them the hyperplane $H$ is
simply tangent to $C$ (since locally $H$ cuts $2p_i$ on $C$). Now
one can just compute (following Andreotti) the Gauss map on the
theta divisor $\Theta[\xi]$ of $J$ around the origin by
parametrizing $C$ around each of the points $p_i$.  In
appropriate coordinates the Gauss map is the gradient of a
function whose leading term cuts out $V_\xi$ (see the proof of
2.3.4). Cf\. \cite{MSV1, p\. 735 (ii)}, \cite{AMSV1, p\. 21,
3.6(iv)}, \cite{AMSV2, tables 1, 2 ($P_8$)}. Q.E.D. \enddemo

In fact, instead of just citing our previous calculations, we
prove the following more general result. Let $A$ be a
$g$-dimensional abelian variety.  A hypersurface germ $(M, 0)
\subset A$ is of {\it translation type\/} if there exist $g - 1$
germs $\Gamma_1$, $\dots$, $\Gamma_{g-1}$ of smooth curves
through $0$ such that the addition map $\mu: \Gamma_1 \times
\dots \times \Gamma_{g-1} \rightarrow A$ has image $(M, 0)$. Such
a germ is of {\it nondegenerate translation type\/} if the curves
$\Gamma_1$, $\dots$, $\Gamma_{g-1}$ can be chosen so that their
tangent lines at $0$ are linearly independent. The theta divisor
$\Theta[\xi]\subset J(C)$ is of nondegenerate translation type at
$0$ for all odd theta characteristics $\xi$ of a generic curve
$C$ of genus $g\ge 3$. (For discussion of theta divisors' being
(singly or doubly) of translation type, cf\. \cite{L}, \cite{M2,
pp\. 81-85}.)

\proclaim{2.3.4 Proposition}  For $(A, \Theta, \xi) \in
\tilde{\Cal A}_g$, $g \ge 4$, if the theta divisor $\Theta[\xi]
\subset A$ is of nondegenerate translation type at $0$, then
$V_\xi \subset \Bbb E_{\xi} \cong \Bbb P^{g-2}$ is a Fermat cubic
hypersurface or a degeneration thereof (including the possibility
that it is undefined, i.e., that $\text{mult}_0\Theta[\xi] > 3$).
\endproclaim

\demo{Proof}  Working locally in $\Bbb C^g$ at $0$, let $\sigma_i
:(\triangle, 0) \rightarrow (\Gamma_i, 0) \subset \Bbb C^g$ be a
local parametrization of the $i\th$ curve germ (by a disc), $i =
1, \dots, g - 1$. Then by a linear change of coordinates we may
arrange that $\sigma_i(t_i)$ has $t_i$ in the $i\th$ coordinate
position and higher order terms in the other coordinates.  Let
$f_i$ denote the last component function (i.e., the $g\th$
coordinate) of $\sigma_i$.  Now if $F(z_1, \dots, z_g)=0$ is a
local analytic equation for $\Theta[\xi]$, then we may uniquely
write $F(z_1, \dots, z_g) = z_g - f(z_1, \dots, z_{g-1})$, where
$f(0) = 0$ and $f$ has no linear term.  Therefore, $f(z_1, \dots,
z_{g-1}) = f_1(t_1) + \dots + f_{g-1}(t_{g-1}) = f_1(z_1) + \dots
+ f_{g-1}(z_{g-1}) + (\text{higher order terms in } z_1,\dots,
z_{g-1})$.  (The Gauss map of $\Theta[\xi] \subset (\Bbb C^g, 0)$
has the local form $\gamma(z_1, \dots, z_{g-1}) = \nabla f(z_1,
\dots, z_{g-1})$.) Since there exists an equation $F$ for
$\Theta[\xi]$ at $0$ that is odd, $f(z_1, \dots, z_{g-1})$ is an
odd function of $z_1, \dots, z_{g-1}$, with no linear term.
Therefore the Taylor expansion of $F$ at $0$ is $$ F(z_1, \dots,
z_g) = z_g - (\lambda_1 z_1^3 + \dots + \lambda_{g-1} z_{g-1}^3)
+ \dots, $$ so $\overline{m} = -(\lambda_1 z_1^3 + \dots
+\lambda_{g-1} z_{g-1}^3)$. Q.E.D. \enddemo

\head 3. Classical invariant theory \endhead

\subhead 3.1\quad Globalization of invariants \endsubhead
\vskip.5\baselineskip

We review the classical invariant theory of degree $k$ forms in
$n$ variables and we show how to express a classical invariant
globally, so that it applies to degree $k$ forms on a vector
bundle. Let $E$ be an $n$-dimensional vector space over $\Bbb C$.
Then the general linear group $\text{GL}(E)$ acts on $F =
S^kE^*$, the vector space of degree $k$ forms (homogeneous
polynomials) on $E$, and hence on $\Cal P = S^d F^*$, the forms
of degree $d$ on $F$.  (For $f \in F$, $x \in E$, $(g \cdot f)(x)
= f(g^{-1}x)$; for $\varphi \in \Cal P$, $f \in F$, $(g \cdot
\varphi)(f) = \varphi(g^{-1} \cdot f)$.)  A degree $d$ {\it
invariant\/} of the degree $k$ forms on $E$ is an element
$\varphi \in \Cal P$ such that $g \cdot \varphi = \varphi$ for
all $g \in \text{SL}(E)$, the special linear group.  It is easy
to determine the action of a general element of $\text{GL}(E)$ on
such $\varphi$.  Namely, for $g \in \text{GL}(E)$, write $g =
(cI)h$, where $c$ is a scalar, $I$ is the identity and $h \in
\text{SL}(E)$; of course, $\det(g) = c^n$.  Then $g \cdot \varphi
= c^{kd}\varphi$.  In particular, $\text{GL}(E)$ acts on the
vector space $\Bbb C \varphi$, so there exists a character $\chi$
of $\text{GL}(E)$ such that, for all $g \in \text{GL}(E)$, $g
\cdot \varphi = \chi(g) \varphi$.  Any character $\chi$ of
$\text{GL}(E)$ (as an algebraic group over $\Bbb C$) has the form
$\chi(g) = \det(g)^w$ for some $w \in \Bbb Z$. Therefore, if the
invariant $\varphi$ is nonzero, then there exists a unique
integer $w$ such that, whenever $g = (cI)h$ as above, we have
$c^{kd} = \det(g)^w$.  It follows that $n|kd$.  Thus, in order
for nontrivial degree $d$ invariants of degree $k$ forms to
exist, we must have $kd = nw$ for a nonnegative integer $w$ (cf\.
\cite{G-Y, Ch. II, \S 31, p\. 28, Ch. XII, \S 199, p\. 246}).

Let $kd = nw$ for a nonnegative integer $w$ and let $$ \Cal
P_{(w)} = \{ \varphi \in \Cal P\ |\ g\cdot \varphi = \det(g)^w
\varphi \text{ for all } g \in \text{GL}(E)\} $$ be the set of
{\it relative invariants\/} of weight $w$.  To create absolute
invariants from relative ones, we must eliminate the determinant
factor.  Thus, we let $W = (\Lambda^n E)^w$ and consider $\Cal P'
= \Cal P \otimes W^*$ with the natural action of $\text{GL}(E)$.
Now $\Cal P_{(w)} \otimes W^* \subset \Cal P'$ consists of
elements that transform trivially under $\text{GL}(E)$, and we
can regard an element $\varphi \in \Cal P'$ as a homogeneous
polynomial of degree $d$ on $F = S^kE^*$ with values in $W^* =
(\Lambda^n E^*)^w$.

Now suppose that our $n$-dimensional vector space $E$ is $\Bbb
C^n$ (i.e., that we have a distinguished basis for $E$).  Then
$\Lambda^n\Bbb C^n$ and hence $W = (\Lambda^n\Bbb C^n)^w$ are
canonically trivial (as $\Bbb C$-vector spaces, not as
$\text{GL}_n(\Bbb C)$-modules).  Let $R_d$ denote the space of
degree $d$ invariants of degree $k$ forms in $n$ variables, $$
R_d = \{\varphi \in S^d(S^k(\Bbb C^{n*})^*)\ |\ g \cdot \varphi =
\varphi \text{ for all } g \in \text{SL}_n(\Bbb C)\}. $$

\proclaim{3.1.1 Proposition}  Let $kd = nw$ for a nonnegative
integer $w$.  For any rank $n$ complex vector bundle $\Cal E$,
there exists a canonical injection $G$ from $R_d$ to the set of
degree $d$ homogeneous polynomial vector bundle mappings from
$S^k\Cal E^*$ to $(\Lambda^n\Cal E^*)^w$. \endproclaim

\demo{Proof}  Let $U$ be an open subset of the base space over
which $\Cal E$ is trivial and choose a trivializing frame
$\sigmabar = \{\sigma_1, \dots, \sigma_n\}$ for $\Cal E^*$ over
$U$.  Consider a section $m = \sum a_{k_1\dots k_n}
\sigma^{k_1}_1 \dots, \sigma^{k_n}_n$ of $S^k\Cal E^*$ over $U$,
where $a_{k_1\dots k_n}$ are regular functions on $U$; let $M =
\sum a_{k_1\dots k_n} x^{k_1}_1 \dots x^{k_n}_n$ be the
corresponding polynomial of degree $k$ (with variable
coefficients).  Then set $G_{\varphi}(m) = \varphi(M)(\sigma_1
\wedge \dots \wedge \sigma_n)^{\otimes w}$.  Now consider another
choice $\tilde{\sigmabar} = \{\tilde{\sigma}_1, \dots,
\tilde{\sigma}_n\}$ of frame for $\Cal E^*$ over $U$. Treating
$\sigmabar$ and $\tilde{\sigmabar}$ as $n$-entry column vectors,
we express $\tilde{\sigmabar} = g \sigmabar$, where $g =
(g_{ij}): U \rightarrow \text{GL}(n, \Bbb C)$ is an invertible
matrix of regular functions.  Then we use the relative invariance
of $\varphi$ to compute that $\varphi(\tilde{M})(\tilde{\sigma}_1
\wedge \dots \wedge \tilde{\sigma}_n)^{\otimes w} =
\varphi(M)(\sigma_1 \wedge \dots \wedge \sigma_n)^{\otimes w}$,
so there is everywhere a unique local realization of
$G_{\varphi}$, hence a unique global realization. Q.E.D.
\enddemo

$G$ is actually a ring homomorphism.  Let $R = \oplus R_d$ be the
graded ring of invariants. Then a degree $d$ element $\varphi \in
R_d$ gives a polynomial vector bundle map $G_{\varphi}: S^k\Cal
E^* \rightarrow (\Lambda^n\Cal E^*)^w$, where $kd =nw$.  The
function $G$ is multiplicative in the sense that for $\varphi \in
R_d$ and $\psi \in R_{d'}$, $G_{\varphi\psi} = G_{\varphi} \cdot
G_{\psi}$ as maps from $S^k\Cal E^*$ to $(\Lambda^n\Cal E^*)^{(w
+ w')}$. The product $G_{\varphi} \cdot G_{\psi}$ is defined by
$(G_{\varphi} \cdot G_{\psi})(m) = G_{\varphi}(m) \cdot
G_{\psi}(m)$, the natural multiplication of sections being
obtained from the tensor product of line bundles $$
(\Lambda^n\Cal E^*)^w \otimes (\Lambda^n\Cal E^*)^{w'}
\overset\sim\to\rightarrow (\Lambda^n\Cal E^*)^{(w + w')}. $$ In
terms of a choice of local frame $\sigmabar$, $$ \align
G_{\varphi\psi}(m) &= (\varphi\psi)(M)(\sigma_1 \wedge
\dots\wedge \sigma_n)^{\otimes(w + w')}\\
&=(\varphi(M)\psi(M))\big((\sigma_1\wedge \dots \wedge
\sigma_n)^{\otimes w} \cdot (\sigma_1 \wedge \dots \wedge
\sigma_n)^{\otimes w'}\big)\\ &=\varphi(M)(\sigma_1
\wedge\dots\wedge \sigma_n)^{\otimes w} \cdot \psi(M)(\sigma_1
\wedge \dots \wedge \sigma_n)^{\otimes w'}\\ &= G_{\varphi}(m)
\cdot G_{\psi}(m). \endalign $$

\example{3.1.2\  Examples}  (i)  The {\it discriminant\/}
$\delta$ is an invariant of degree $n(k-1)^{(n-1)}$, as a
polynomial in the coefficients of degree $k$ forms in $n$
variables. It is defined by $\delta(f) = R(\frac{\partial
f}{\partial x_1}, \dots, \frac{\partial f}{\partial x_n})$, where
$R$ is the $n$-variable resultant \cite{vW, p\. 15}, and it has
the property that $\delta(f) = 0$ if and only if the projective
hypersurface $f = 0$ is singular (as a subscheme of $\Bbb
P^{n-1}$) (cf\. \cite{M-F, Prop\. 4.2, p\. 79}).

\noindent(ii)  All invariants of positive degree vanish on forms
which are missing a variable. That is, suppose that $f$ is a
degree $k$ form in $n$ variables such that, after some linear
change of coordinates, $f$ can be expressed in terms of fewer
variables; geometrically, the projective hypersurface $f = 0$ in
$\Bbb P^{n-1}$ is a cone. Then any invariant $\varphi \in R_d$,
$d > 0$, must vanish on $f$; i.e., such an $f$ is a {\it
nullform\/} (or is {\it unstable\/} (cf\. \cite{M-F, Thm\. 2.1,
p\. 49}).  Indeed, assuming that the variable $x_n$ is missing,
consider the 1-parameter group $\lambda: \Bbb C^* \rightarrow
\text{SL}(S^k\Bbb C^{n*})$ induced on degree $k$ forms by the
action $t \cdot (x_1, \dots, x_{n-1}, x_n) = (tx_1, \dots,
tx_{n-1}, t^{-(n-1)}x_n)$, $t \in \Bbb C^*$.  Then $\lambda(t)(f)
= t^kf$, so $\varphi(f) = \varphi(\lambda(t)(f)) =
t^{kd}\varphi(f)$.  Hence, letting $t \rightarrow 0$, we get
$\varphi(f) = 0$. \endexample

\vskip.5\baselineskip \subhead 3.2\quad Invariants of cubic forms
\endsubhead \vskip.5\baselineskip

In our applications we will use the invariant theory of cubic
forms in $g-1$ variables; thus $k = 3$ and $n = g - 1$, so the
fundamental equation becomes $3d = (g - 1)w$.  Now $E = \Bbb
C^{g-1}$ and $\text{GL}(E)$ acts on $F = S^3E^*$, the vector
space of cubic forms on $E$, and hence on $S^dF^*$, the
homogeneous polynomials of degree $d$ on $F$. In this situation,
the discriminant $\delta$ has degree $(g - 1)2^{g - 2}$.

The {\it standard Fermat cubic form\/} in $g-1$ variables is $f_0
= x_1^3 + \dots + x_{g-1}^3 \in F$ and the zero set of $f_0$ in
$\Bbb P^{g-2}$ is the standard Fermat cubic hypersurface.  The
$\text{GL}_{g-1}(\Bbb C)$-orbit of $f_0$ consists of all cubic
forms that can be expressed as the sum of cubes of $g - 1$
linearly independent linear forms, and we will refer to any
element of this orbit as a {\it Fermat cubic\/} in $g-1$
variables.  Let $I \subset R$ be the homogeneous ideal of all
invariants vanishing on $f_0$ (or equivalently on the orbit
$\text{GL}(E) \cdot f_0$, or on the closure of this orbit).

We summarize the well known structure of the ring of invariants
of ternary cubics and the (principal) ideal of all invariants
vanishing on the Fermat cubic $x^3 + y^3 + z^3$.  For $E = \Bbb
C^3$, $F = S^3E^*$ is $10$-dimensional, and we are interested in
the graded ring $R$ of $\text{SL}_3(\Bbb C)$-invariant polynomial
functions on $F$.  The structure theorem for invariants of
ternary cubics states that there are two algebraically
independent invariants  $S$ and $T$, of degrees $4$ and $6$
respectively, such that $R = \Bbb C[S, T]$ (\cite{E1,
\S\S291-293, 295, pp\. 381-389}, \cite{J, p\. 4}).

We now describe the explicit forms for the discriminant $\delta$
and the ideal $I$ of invariants vanishing on the Fermat cubic.
Using Hesse's canonical form for cubics (cf\. \cite{E1, \S229},
\cite{J, p\. 6}), any nonsingular cubic in $x$, $y$ and $z$ can
be put, by a linear change of coordinates, in the form $x^3 + y^3
+ z^3 + 6mxyz$; in particular, any invariant of cubics is
determined by its values on such normalized cubics.  Now $S = m -
m^4$ and $T = 1 - 20m^3 - 8m^6$ \cite{E1, pp\. 384, 386},
\cite{J, p\. 7}. (For further information, see \cite{D-K, p\.
250}, \cite{J, p\. 6}, \cite{E1, p\. 384}.)  The discriminant
$\delta$ is $T^2 + 64S^3$ \cite{J, p\. 27} and the ideal $I$ is
$(S)$ (cf\. \cite{D-K, p\. 251}, \cite{G-Y, \S 248, pp\.
312-313}).  All members of the closure of the Fermat orbit
$\text{GL}_3(\Bbb C) \cdot (x^3 + y^3 + z^3)$ are described
geometrically in \cite{D-K, Prop\. 5.13.2 (ii), p\. 251}.
Finally, the $j$-invariant, the fundamental rational invariant of
cubic curves in $\Bbb P^2$, is given by $j =
\text{(constant)}\,S^3/\delta$.

\eject \head 4. The construction of certain Siegel modular forms
\endhead

\subhead 4.1\quad The ring homomorphism $h_{\xi}$ \endsubhead
\vskip.5\baselineskip

We assemble the results of the previous sections to construct a
ring homomorphism from invariants of cubic forms to Siegel
modular forms. Let $\overline{\Bbb N} = \{0, 1, 2, \dots \}$, the
additive semigroup of (extended) natural numbers. Let $R = \oplus
R_d$, $d \in \overline{\Bbb N}$, be the graded ring of invariants
of cubic forms in $g-1$ variables.  Let $\Cal S = \oplus \Cal
S_k$, $k \in \frac12\overline{\Bbb N}$, be the graded ring of
genus $g$ Siegel modular forms with respect to $\Gamma(4, 8)$.

A {\it Siegel modular form\/} of genus $g > 1$, weight $k \in
\frac12\overline{\Bbb N}$ and level $\Gamma(4, 8)$ (with trivial
character) is a holomorphic function $f: \Cal H_g \rightarrow
\Bbb C$ such that for every $\gamma = \left(\smallmatrix A & B\\
C & D\endsmallmatrix\right) \in \Gamma(4, 8)$ and $\Omega \in
\Cal H_g$, $f(\gamma \cdot \Omega) = \det(C\Omega +
D)^kf(\Omega)$.  Thus, for $k \in \frac12\overline{\Bbb N}$,
$\Cal S_k = H^0(\Cal A'_g, P^{2k})$.

\proclaim{4.1.1 Theorem}  For every genus $g \ge 4$ and each
choice of odd theta characteristic $\xi$, there exists a ring
homomorphism $$ h_{\xi}: R \rightarrow \Cal S $$ with the
following properties \roster \item"{(0)}" $h_{\xi}(R_d) \subset
\Cal S_{md}$, where $m = \frac{g+8}{2(g-1)}$. \item"{(1)}"  If
$\delta \in R$ is the discriminant for cubic forms, then
$h_{\xi}(\delta) \ne 0$. \item"{(2)}"  If $I \subset R$ is the
homogeneous ideal of all invariants vanishing on the Fermat
cubic, then all elements of $h_{\xi}(I) \subset \Cal S$ vanish on
period matrices of genus $g$ Riemann surfaces. \endroster
\endproclaim

\remark{\rm{4.1.2}\ Remarks} (i)  Since $3d = (g -1)w$ for some
$w \in \overline{\Bbb N}$ in order for nontrivial invariants to
exist, it follows that if $R_d \ne 0$, then $md \in
\frac12\overline{\Bbb N}$.

\noindent (ii)  $h_{\xi}(\delta) \in \Cal S_k$ for $k = (g +
8)2^{g-3}$. \endremark

\demo{Proof}  From (2.1) we have an open subset $\Cal U = \Cal
U_{\xi}$ of the moduli space $\Cal A'_g$, and over $\Cal U$ a
rank $g-1$ vector bundle $\Cal E = \Cal E_{\xi} \subset \Cal V$,
where $\Cal V = H^*$ is the dual of the Hodge vector bundle $H$.
In (2.2) we constructed a homomorphism $\overline{m}_{\xi}: N
\rightarrow S^3\Cal E^*$ from the fundamental negative line
bundle $N$ to the bundle of cubic forms on $\Cal E$.  Now let
$\varphi$ be a degree $d$ invariant of cubic forms in $g-1$
variables, i.e., $\varphi \in R_d = S^d((S^3\Bbb
C^{g-1*})^*)_{(w)}$  where $w = 3d/(g - 1)$.  By (3.1.1),
$\varphi$ determines a degree $d$ homogeneous polynomial mapping
$G_{\varphi}: S^3\Cal E^* \rightarrow (\Lambda^{g-1}\Cal E^*)^w$.
Therefore we obtain (on $\Cal U$) a degree $d$ homogeneous
polynomial mapping $G_{\varphi} \comp \overline{m}_{\xi} : N
\rightarrow (\Lambda^{g-1}\Cal E^*)^w$, hence a (linear)
homomorphism from $N^d$ to $(\Lambda^{g-1}\Cal E^*)^w$, and so a
section $s_{\varphi}$ of the line bundle $N^{-d} \otimes
(\Lambda^{g-1}\Cal E^*)^w$ on $\Cal U$.  That is, $\varphi \in
R_d$ determines $s_{\varphi} \in H^0(\Cal U, N^{-d} \otimes
(\Lambda^{g-1}\Cal E^*)^w)$.

{}From the exact sequence $0 \rightarrow \Cal L_{\xi} \rightarrow H
\rightarrow \Cal E_{\xi}^* \rightarrow 0$ on $\Cal U$, we have
$\Lambda^g H \cong \Cal L_{\xi} \otimes \Lambda^{g-1}\Cal
E^*_{\xi}$. Thus $\Lambda^{g-1}\Cal E^*_{\xi} \cong \Cal
L_{\xi}^{-1} \otimes \Lambda^g H \cong N^{-1} \otimes \Lambda^g H
\cong P^3$.  Thus $N^{-d} \otimes(\Lambda^{g-1}\Cal E^*)^w \cong
P^{(d+3w)}$.  Since $w = 3d/(g-1)$, we have $d + 3w = \mu d \in
\overline{\Bbb N}$, where $\mu = (g + 8)/(g - 1)$.  Thus, from a
degree $d$ invariant $\varphi$ we get a section $s_{\varphi}$ of
$N^{-\mu d} = P^{\mu d}$ on $\Cal U$.  By (2.1.3), $H^0(\Cal U,
P^{\mu d}) \cong H^0(\Cal A'_g, P^{\mu d})$, and we let
$h_{\xi}(\varphi)$ denote the unique extension of $s_{\varphi}$
to a section of $P^{\mu d}$ over all of $\Cal A'_g$. Therefore
$h_{\xi}(\varphi)$ is a Siegel modular form of weight $k = \mu
d/2 = md$, where $m$ is as stated in (0).

Next we check that $h_{\xi}(\varphi \psi) =
h_{\xi}(\varphi)h_{\xi}(\psi)$ for $\varphi \in R_d$, $\psi \in
R_{d'}$.  Let $\tilde{G}_{\varphi}$ denote the homomorphism
$S^d(S^3\Cal E^*) \rightarrow (\Lambda^{g-1}\Cal E^*)^w$
corresponding to the degree $d$ homogeneous polynomial map
$G_{\varphi}: S^3\Cal E^* \rightarrow (\Lambda^{g-1}\Cal E^*)^w$
defined by $\varphi$.  Then $h_{\xi}(\varphi)$ is the section of
$N^{-2md} \cong N^{-d} \otimes (\Lambda^{g-1}\Cal E^*)^w$ which
corresponds to the composite homomorphism $\tilde{G}_{\varphi}
\comp(\overline{m}_{\xi})^d: N^d \rightarrow S^d(S^3\Cal E^*)
\rightarrow (\Lambda^{g-1}\Cal E^*)^w$. Thus the equality of the
product $h_{\xi}(\varphi)h_{\xi}(\psi)$ in\break $(N^{-d} \otimes
(\Lambda^{g-1}\Cal E^*)^w) \otimes (N^{-d'} \otimes
(\Lambda^{g-1}\Cal E^*)^{w'})$ with $h_{\xi}(\varphi \psi)$ in
$N^{-(d+d')} \otimes (\Lambda^{g-1}\Cal E^*)^{(w + w')}$ follows
from the equality of the product homomorphism $$
(\tilde{G}_{\varphi} \comp
(\overline{m}_{\xi})^d)(\tilde{G}_{\psi}
\comp(\overline{m}_{\xi})^{d'}): N^d \otimes N^{d'} \rightarrow
(\Lambda^{g-1}\Cal E^*)^w \otimes (\Lambda^{g-1}\Cal E^*)^{w'} $$
with $\tilde{G}_{\varphi \psi} \comp (\overline{m})_{\xi}^{(d +
d')}$.  But $(\Tilde{G}_{\varphi} \comp
(\overline{m}_{\xi})^d)(\tilde{G}_{\psi}
\comp(\overline{m}_{\xi})^{d'}) = (\tilde{G}_{\varphi} \cdot
\tilde{G}_{\psi}) \comp (\overline{m}_{\xi})^{(d + d')}$, so it
suffices to check that the product $$ G_{\varphi} \cdot
\tilde{G}_{\psi}: S^{(d+d')}(S^3\Cal E^*) \rightarrow
(\Lambda^{g-1}\Cal E^*)^{(w + w')} $$ agrees with
$\tilde{G}_{\varphi \psi}$.  But $\tilde{G}_{\varphi} \cdot
\tilde{G}_{\psi}$ and $\tilde{G}_{\varphi \psi}$ agree since they
are the homomorphisms corresponding to the degree $d + d'$
homogeneous polynomial maps $G_{\varphi} \cdot G_{\psi}$ and
$G_{\varphi \psi} : S^3\Cal E^* \rightarrow (\Lambda^{g-1}\Cal
E^*)^{(w + w')}$, respectively, and we know from \S 3 that
$G_{\varphi} \cdot G_{\psi} = G_{\varphi \psi}$.

To establish properties (1) and (2), it suffices to recall from
(2.3.2) (ii) that for a generic genus $g$ Riemann surface and any
choice of odd theta characteristic $\xi$, the associated cubic
form $\overline{m}_{\xi}(z)$ is a Fermat cubic in $g - 1$
variables.  Indeed, if $\Omega_0 \in \Cal H_g$ is a period matrix
over a generic Jacobian $(J, \Theta, \xi') \in \Cal A'_g$, then
(1) $(h_{\xi}(\delta))(\Omega_0) = \delta(\overline{m}_{\xi}(z))
\ne 0$ since the Fermat cubic hypersurface is nonsingular, and
(2) holds by construction. Q.E.D. \enddemo

\remark{\rm{4.1.3}\ Remark}  The construction of the graded ring
homomorphism $h_{\xi}$ and the calculation of its degree is
simplified slightly by the introduction of a formal cube root
$P^{1/3}$ of $P$ (in analogy with the squaring principle of
\cite{H-T, p\. 77}).  Then $P \otimes S^3\Cal E^* = S^3(P^{1/3}
\otimes \Cal E^*)$, so $G_{\varphi}$ takes values in
$\Lambda^{g-1}(P^{1/3} \otimes \Cal E^*)^w \cong P^{(d+3w)}$.
\endremark

\vskip.5\baselineskip \subhead 4.1.4 \ \  Reprise \endsubhead
Given $\Omega \in \Cal H_g$, let us try to actually get a number
from the above abstract construction.  Recall that we have a
linear form $\ell$  and a cubic form $m$ on the standard vector
space $\Bbb C^g$.  Let $\varphi$ be a classical invariant of
cubic forms in $g - 1$ variables. Assuming the linear form $\ell$
is not identically $0$, let $E \subset \Bbb C^g$ be the
hyperplane $\ell(z) = 0$ and let $\overline{m}$ be the
restriction of the cubic form to $E$.  Then if we take any basis
$B = \{v_2, \dots, v_g\}$ for $E$, we get a cubic form $M_B(x) =
\overline{m}(x_2v_2 + \dots +x_gv_g)$ in $g - 1$ variables, so we
can apply $\varphi$ to get a number $\varphi(M_B)$.  If we change
the basis of $E$ by $\alpha \in \text{GL}(E)$, then
$\varphi(M_{\alpha B}) = \det(\alpha)^p \cdot (M_B)$ for some
fixed integer $p$.  Now suppose we construct the basis for $E$
only in the following way.  We take a basis $\Cal B$ for $\Bbb
C^{g*}$ of the form $\{\ell_1 = \ell, \ell_2, \dots, \ell_g\}$.
Then the dual basis $v_1, v_2, \dots, v_g$ for $\Bbb C^g$ has the
property that $v_2, \dots, v_g$ form a basis $B$ for $E$ and we
consider, instead of $\varphi(M_B)$, the number $\det(\Cal B)^p
\cdot \varphi(M_B)$.  It is easy to compute that if we change to
a new basis $\tilde{\Cal B} = \{\ell_1 = \ell, \tilde{\ell}_2,
\dots, \tilde{\ell}_g\}$ then $\det(\tilde{\Cal B})^p \cdot
\varphi(M_{\tilde B}) = \det(\Cal B)^p \cdot \varphi(M_B)$.  Thus
the numerical value $\det(\Cal B)^p \cdot \varphi(M_B)$ does not
depend on the choice of (special) basis $\Cal B$ of $\Bbb
C^{g*}$; hence this formula can be used to evaluate the modular
form $h_{\xi}(\varphi)$ at a period matrix $\Omega$ directly in
terms of the Taylor expansion of the theta function about $z =
0$.

\vskip.5\baselineskip \subhead 4.2\quad The genus $4$ case
\endsubhead \vskip.5\baselineskip

Let $\script J'_4$ be the preimage in $\Cal A'_4$ of the locus of
genus $4$ Jacobians $\script J_4 \subset \Cal A_4$ and let
$\bar{\script J}'_4$ be the closure of $\script J'_4$ in $\Cal
A'_4$.

\proclaim{4.2.1 Proposition}  For $g = 4$, let $h$ be the
homomorphism determined by a marked odd theta characteristic, let
$S$ be the classical quartic invariant of ternary cubics, and set
$f = h(S)$.  Then $f$ is a Siegel modular form of weight $8$ with
respect to $\Gamma' = \Gamma(4, 8)$ and $\bar{\script J}'_4
\subset \Cal A'_4$ is the zero divisor of $f$ (with multiplicity
1). \endproclaim

\demo{Proof}  When $g = 4$ we have $m = 2$ in (4.1.1), and hence
if $\varphi = S$ is the classical quartic invariant of cubic
forms on $\Bbb C^3$, then $f = h(S)$ is a Siegel modular form of
degree $8$. From our previous indirect argument \cite{AMSV2, pp\.
5-6}, we know that this modular form $f$ is not identically zero.
(We also have \cite{MSV2} a direct local argument that the local
modulus $a(\Omega)$ in the cubic form $\overline{m} = x^3 + y^3 +
z^3 + a(\Omega)xyz$, is not identically $0$.  Presumably it could
also be verified by numerical calculation that $f(\Omega)$ is
nonzero for a specific period matrix $\Omega$.)  Therefore, it
remains only to use the following. \enddemo

\proclaim{4.2.2 Lemma}  The divisor class of $\bar{\script J}'_4$
in $\text{Pic}(\Cal A'_4) \otimes \Bbb Q$ is equal to $8\lambda$,
where $\lambda = c_1(\Lambda)$ is the class of the Hodge line
bundle. \endproclaim

\demo{Proof \rm{(Mumford)}} We use adjunction and Mumford's
formula for the canonical bundle of $\Cal M_g$.  Note that $\Cal
A'_4$ is smooth, $\script J'_4$ is a smooth
Zariski-locally-closed divisor of $\Cal A'_4$, and $\script J'_4$
is closed in  $\Cal A'_4 - \Cal R'$, where the subset $\Cal R'$
of decomposables is Zariski-closed in $\Cal A'_4$ of codimension
$\ge 2$ (and is actually contained in $\bar{\script J}'_4$).  Now
we consider the conormal bundle sequence for $\script J = \script
J'_4 \subset \Cal A'_4 - \Cal R' = \Cal A:$ $$ 0 \rightarrow
N^*(\script J/\Cal A) \rightarrow T^*(\Cal A)|\script J
\rightarrow T^*(\script J) \rightarrow 0. $$ Hence, $K_{\Cal
A|\script J} = K_{\script J} + [N^*]$. By the standard formula
$K_{\Cal A'_g} = (g + 1) \lambda$ (cf\. \cite{Fr 2. III \S 2}),
we have $K_{\Cal A} = 5 \lambda$, and by Mumford's formula
$K_{\Cal M_g} = 13 \lambda$ (cf\. \cite{H-M, \S 2}), we have
$K_{\script J} = 13\overline{\lambda}$, where
$\overline{\lambda}$ is the restriction of $\lambda$ to $\script
J$. Therefore $[N^*] = 5 \overline{\lambda} -
13\overline{\lambda} = -8\overline{\lambda}$, so $[N] =
8\overline{\lambda}$.  By Borel's results \cite{Bo}, $\Cal A'_g$
has Picard number one for $g \ge 4$, so $[\script J] = r\lambda$
in $\text{Pic}(\Cal A) \otimes \Bbb Q$ for some $r \in \Bbb Q$.
Thus, if $i:\script J \hookrightarrow \Cal A$ is the inclusion,
we get $i^*[\script J] = i^*(r\lambda) = r\overline{\lambda}$ in
$\text{Pic}(\script J) \otimes \Bbb Q$.  But of course
$i^*[\script J] = i^*[\Cal O_{\Cal A}(\script J)] = [i^*\Cal
O_{\Cal A}(\script J)] = [N]$.  Thus we have both $[N] =
8\overline{\lambda}$ and $r\overline{\lambda} = [N]$, whence $r =
8$.  Therefore, $[\bar{\script J}'_4] = 8 \lambda$ in
$\text{Pic}(\Cal A'_4) \otimes \Bbb Q$.  Q.E.D. (4.2.2)

We can conclude the genus $4$ argument by again using Picard
number one.  Let $D$ be the divisor of $f$ in $\Cal A'_4$.  Since
$f$ vanishes on $\script J'_4$ and hence on $\bar{\script J}'_4$,
we have $D = \bar{\script J}'_4 + E$ with $E \ge 0$.  Now take
classes in $\text{Pic} \otimes \Bbb Q$.  We have $8\lambda = [D]
= [\bar{\script J}'_4] + [E] = 8\lambda + [E]$, so $[E] = 0$.
But the divisor class (in $\text{Pic}$) of a nonzero effective
divisor on $\Cal A'_4$ cannot be $0$ (by the existence of the
Satake compactification $(\Cal A'_4)^*$, which has small
boundary), and hence a nonzero effective divisor on $\Cal A'_4$
cannot be torsion either.  Therefore the divisor $E$ is $0$ and
hence $D = \bar{\script J}'_4$, as desired.  Q.E.D. (4.2.1)
\enddemo

We can reformulate the previous result in the following way.

\proclaim{4.2.3 Corollary}  Let $j: \tilde{\Cal A}_4
\dashrightarrow \Bbb P^1$ be the rational function defined by
$j(A, \Theta, \xi) = j(a)$, where (after a suitable linear change
of coordinates, if possible) $\overline{m}_{\xi}(x, y, z) = x^3 +
y^3 + z^3 + a \cdot xyz$ for some $a \in \Bbb C$, and $j(a)$ is
the classical $j$-invariant of $\overline{m}_{\xi}$.  The zero
divisor of $j$ in $\tilde{\Cal A}_4$ is the closure of
$\tilde{\script J}_4$ (with multiplicity $3$). \endproclaim

\proclaim{4.2.4 Corollary}  The modular form $f$ of (4.2.1) is a
nonzero scalar multiple of Schottky's equation.  A fortiori, $f$
has level one; i.e., $f$ is a modular form of weight $8$ with
respect to the entire group $\text{\rm Sp}(2g, \Bbb Z)$.
\endproclaim

\demo{Proof}  Let $J$ be the Schottky equation, which is a genus
$4$ Siegel modular form of weight $8$ and level $1$, and which
has zero divisor $\bar{\script J}_4 \subset \Cal A_4$ (with
multiplicity one) (\cite{I4}, \cite{Fr1}).  In particular, $J$
can be viewed on $\Cal A'_4$ as a section of $P^{16} = \Lambda^8$
with zero divisor $\bar{\script J}'_4$.  Thus, $f$ and $J$ are
sections of the same line bundle and have the same zero divisor
on $\Cal A'_4$, and hence $f/J$ is a nowhere vanishing
holomorphic function $c$ on $\Cal A'_4$.  But  $H^0(\Cal A'_4,
\Cal O) = \Bbb C$ (by the existence of the Satake
compactification $(\Cal A'_4)^*$, with boundary $(\Cal A'_4)^* -
\Cal A'_4$ of codimension $\ge 2$), so $c$ is a nonzero scalar.
Then, since $f = cJ$, we see that in fact $f$ has level $1$.
Q.E.D. \enddemo

Here is a self-contained argument for level one, based directly
on Proposition (4.2.1).

\proclaim{4.2.5 Lemma}  Let $f$ be a modular form on $\Cal H_g$,
$g \ge 3$, of integral weight $k$ and level $\Gamma$, a finite
index subgroup of $\text{\rm Sp}(2g, \Bbb Z)$.  If the zero
divisor of $f$ is invariant under $\text{\rm Sp}(2g, \Bbb Z)$,
then $f$ is modular with respect to all of $\text{\rm Sp}(2g,
\Bbb Z)$. \endproclaim

\demo{Proof}  Define an action of $\text{Sp}(2g, \Bbb Z)$ on
$\Cal O(\Cal H_g)$ as follows (cf\. \cite{Fr2, p\. 53}).  For
each integer $k$, $f\in\Cal O(\Cal H_g)$ and $\gamma =
\left(\smallmatrix A & B \\ C & D \endsmallmatrix\right) \in
\text{Sp}(2g, \Bbb Z)$, define $(f|_k\gamma)(\Omega) =
\lambda(\gamma, \Omega)^{-k} f(\gamma \cdot \Omega)$, where
$\lambda(\gamma, \Omega) = \det(C \Omega + D)$ and $\gamma \cdot
\Omega = (A \Omega + B)(C \Omega + D)^{-1}$.  If the zero divisor
of $f$ is invariant under $\gamma \in \text{Sp}(2g, \Bbb Z)$,
there is a nowhere zero holomorphic function $u_{\gamma}$ on
$\Cal H_g$ such that $f|_k\gamma = u_{\gamma} \cdot f$.  In this
way, from $f$ we get a 1-cocycle $\{u_{\gamma}\}$ of
$\text{Sp}(2g, \Bbb Z)$ with coefficients in $\Cal O^*(\Cal H_g)$
(a multiplicative group on which $\text{Sp}(2g, \Bbb Z)$ acts
nontrivially).  Now if $f$ is a modular form on $\Cal H_g$ (of
weight $k$) with respect to $\Gamma$, then the functions
$u_{\gamma}$ are holomorphic on $\Cal H_g$ and invariant under
$\Gamma$, and the only modular forms of weight $0$ with respect
to $\Gamma$ are the constants.  Thus $\{u_{\gamma}\} \in
Z^1(\text{Sp}(2g, \Bbb Z), \Bbb C^*) = \text{Hom}(\text{Sp}(2g,
\Bbb Z), \Bbb C^*)$ since $\text{Sp}(2g, \Bbb Z)$ acts trivially
on the coefficients $\Bbb C^*$. But $\text{Hom}(\text{Sp}(2g,
\Bbb Z), \Bbb C^*)$ is trivial, since by \cite{Ma} the commutator
subgroup of $\text{Sp}(2g, \Bbb Z)$ is all of $\text{Sp}(2g, \Bbb
Z)$  if $g \ge 3$. Therefore $u_{\gamma} = 1$ for all $\gamma \in
\text{Sp}(2g, \Bbb Z)$; i.e., $f|_k\gamma = f$ for all $\gamma
\in \text{Sp}(2g, \Bbb Z)$, so $f$ is modular with respect to all
of $\text{Sp}(2g, \Bbb Z)$.  Q.E.D. \enddemo

\vskip.5\baselineskip \subhead 4.3\quad The Fermat locus
\endsubhead \vskip.5\baselineskip

Let $\Cal I_{\xi}$ denote the ideal $(h_{\xi}(I)) \subset \Cal
S$. Every homogeneous element of $\Cal S$, as a holomorphic
function on $\Cal H_g$, has a zero set which is invariant under
$\Gamma' = \Gamma(4, 8)$, so the zero set of the ideal $\Cal
I_{\xi}$ is well defined in $\Cal A'_4 = \Cal H_g/\Gamma'$.  Now
this zero set is actually well defined in $\tilde{\Cal A}_g$,
i.e., modulo $\Gamma[\xi]$.  Note that $\Cal A'_g \rightarrow
\Cal A_g$ is Galois (since $\Gamma'$ is a normal subgroup of
$\text{Sp} (2g, \Bbb Z)$ \cite{I1, Lemma 1, p\. 222}); hence the
finite cover $\pi_{\xi}: \Cal A'_g \rightarrow \tilde{\Cal A}_g$
is also Galois and the Galois group of $\pi_{\xi}$ is
$\Gamma[\xi]/\Gamma'$ (modulo the subgroup $\{\pm I\}$). In
particular, $\Gamma[\xi]$ acts on $\Cal A'_g$ and preserves the
open set $\Cal U\subset\Cal A'_g$.

\proclaim{4.3.1 Lemma}  The zero set in $\Cal H_g$ of each
homogeneous element of $h_{\xi}(I)$ is invariant under the action
of the group $\Gamma[\xi]$, i.e., is well-defined in $\tilde{\Cal
A}_g = \Cal H_g/\Gamma[\xi]$. \endproclaim

\demo{Proof}  Let $I$ be the ideal of all invariants of cubic
forms in $g-1$ variables which vanish on the Fermat cubic. Under
an element of $\Gamma[\xi]$, the theta function
$\vartheta[\xi](z, \Omega)$ transforms to a multiple of itself.
Consequently the action of $\Gamma[\xi]$ preserves the condition
that $I$ vanishes on the cubic $\overline m$ defined by
$\vartheta[\xi](z, \Omega)$. Q.E.D. \enddemo

For each $g \ge 4$, the {\it Fermat locus\/} is the zero locus
$\Cal F_g = V(\Cal I_{\xi}) \subset \tilde{\Cal A}_g$.  If
$\tilde{\script J}_g \subset \tilde{\Cal A}_g$ is the preimage of
the locus of Jacobians $\script J_g \subset \Cal A_g$, then
$\tilde{\script J}_g \subset \Cal F_g$ by part (2) of Theorem
(4.1.1).

\proclaim{4.3.2 Proposition}  Let $\Cal R_g = \cup \Cal A_{ij}
\subset \Cal A_g$ be the locus of all products of (positive
dimensional) principally polarized abelian varieties and let
$\tilde{\Cal R}_g$ be its preimage in $\tilde{\Cal A}_g$.  Then
for each $g \ge 4$, $\tilde{\Cal R}_g \subset \Cal F_g$.
\endproclaim

\demo{Proof}  For a product principally polarized abelian variety
$(A, \Theta) = (A_1, \Theta_1) \times (A_2, \Theta_2)$, every odd
theta characteristic has the form $\xi = \xi_1 + \xi_2$, where
$\xi_1$ and $\xi_2$ are theta characteristics of $(A_1,
\Theta_1)$ and $(A_2, \Theta_2)$, respectively, of opposite
parity. Then $\vartheta[\xi](z) = \vartheta[\xi_1](z_1) \cdot
\vartheta[\xi_2](z_2)$ where $z = z_1 + z_2$, so we can obtain
the Taylor expansion of $\vartheta[\xi](z)$ from the product of
those of $\vartheta[\xi_1](z_1) $ and $\vartheta[\xi_2](z_2)$.
Say that $\xi_1$ is even and $\xi_2$ is odd.  Then $$ \align
\vartheta[\xi](z) &= \vartheta[\xi_1](z_1) \cdot
\vartheta[\xi_2](z_2) \\ &= (c + q(z_1) + \dots )\cdot (\ell(z_2)
+ m(z_2) + \dots )\\ &= c \cdot \ell(z_2) + (c \cdot m(z_2) +
q(z_1) \cdot \ell(z_2)) + \dots\ . \endalign $$ Thus we have
$\ell_{\xi}(z) = c \cdot \ell(z_2)$ and $m_{\xi}(z) = c \cdot
m(z_2) + q(z_1) \cdot \ell(z_2)$.  Hence, assuming $c \ne 0$ and
$\ell \ne 0$ (say that $(A_i, \Theta_i)$ is generic in $\Cal
A_{g_i}$), we see $\overline{m}_{\xi}(z)$ is a scalar multiple of
$m(z_2)$ on the hyperplane $\ell(z_2) = 0$ in $(z_1, z_2)$-space.
In particular, the cubic hypersurface $\overline{m}_{\xi}(z) = 0$
in $\Bbb P(\ell(z_2) = 0)$ is a cone.  But by (3.1.2)(ii), all
(positive degree) invariants of cubic forms in $g-1$ variables
vanish on a cubic cone. Therefore, for the product $(A, \Theta,
\xi) = \big((A_1, \Theta_1) \times (A_2, \Theta_2), \xi_1 +
\xi_2\big) \in \tilde{\Cal A}_g$ with cubic form
$\overline{m}_{\xi}$ and any invariant $\varphi \in R_d$, $d >
0$, we have $\varphi(\overline{m}_{\xi}) = 0$, so $(A, \Theta,
\xi) \in \Cal F_g$. Q.E.D. \enddemo

\proclaim{4.3.3 Corollary}  The ideal $\Cal I_{\xi} =
(h_{\xi}(I)) \subset \Cal S$ consists of cusp forms, i.e.,
elements in the kernel of the Siegel operators for $\Gamma'$
(cf\. \cite{Fr2, Ch. II, \S 6, esp\. Def\. 6.9, p\. 129}).
\endproclaim

\demo{Proof}  Let $(\Cal R'_g)^*$ be the closure of the preimage
of $\Cal R_g$ in the Satake compactification $(\Cal A'_g)^* \cong
\text{Proj}(\Cal S)$ of $\Cal A'_g$.  From the fact that $\Cal
R^*_g$ contains the (single component) Satake boundary $\Cal B =
\Cal A^*_g - \Cal A_g$ of $\Cal A_g$, it follows that $(\Cal
R_g')^*$ contains the Satake boundary $\Cal B' = (\Cal A'_g)^* -
\Cal A'_g$. By (4.3.1) and (4.3.2), we know that $\Cal R'_g
\subset V(\Cal I_{\xi})$ in $\Cal A'_g$, hence $(\Cal R'_g)^*
\subset V(\Cal I_{\xi})$ in $(\Cal A'_g)^*$. Therefore $\Cal B'
\subset V(\Cal I_{\xi})$ in $(\Cal A'_g)^*$; that is, all the
modular forms in $\Cal I_{\xi}$ vanish on the Satake boundary
$\Cal B'$, so they are cusp forms for $\Gamma'$, as claimed.
Q.E.D. \enddemo

\remark{\rm{4.3.4}\ Remark}  In fact, if $R_+=\bigoplus_{d>0}
R_d$, then the ideal $(h_{\xi}(R_+)) \subset \Cal S$ consists of
cusp forms: $h_{\xi}(\delta)$ is a cusp form, by the same
argument as given for (4.3.2), and $R_+$ is the radical of the
ideal $(I, \delta)$ of $R$. \endremark

\remark{\rm{4.3.5}\ Remark}  The boundary of the Satake
compactification $\tilde{\Cal A}^*_g$ of $\tilde{A}_g$ has two
irreducible components ($g \ge 2$), and the closure of $\Cal F_g$
in $\tilde{A}^*_g$ contains both of them by (4.3.3). \endremark

\head 5. Discussion \endhead

If $\pi: \tilde{\Cal A}_g \rightarrow \Cal A_g$ is the map to
level $1$ moduli, we define the {\it big Fermat locus\/} $\Cal
F_g^{\text{(big)}}=\pi(\Cal F_g) \subset \Cal A_g$ and the {\it
small Fermat locus\/} $\Cal F_g^{\text{(small)}} =\{y\ |\
\pi^{-1}(y) \subset \Cal F_g\}\subset \Cal A_g$.  Thus we have
$\script J_g \cup \Cal R_g \subset \Cal F_g^{\text{(small)}}$ for
every $g \ge 4$, and $\bar{\script J}_4 = \Cal
F_4^{\text{(big)}}$.  As a focal point for further work, the
following seems reasonable.

\proclaim{5.1 Conjecture}  For every genus $g \ge 4$,
$\bar{\script J}_g$ is an irreducible component of $\Cal
F_g^{\text{(big)}} \subset \Cal A_g$. \endproclaim

\remark{\rm{5.2}\ Remark}  Let $\Cal C\subset \Cal A_5$ be the
locus of intermediate Jacobians of cubic $3$-folds, and let
$\tilde{\Cal C} \subset \tilde{\Cal A}_5$ denote the canonical
lifting of $\Cal C$ defined by the distinguished odd theta
characteristic of the intermediate Jacobian of a cubic $3$-fold
(cf\. \cite{Do2, p\. 211}). Then it seems to be true, in analogy
with \cite{Do2}, that $\tilde{\Cal C} \subset \Cal F_5$, so $\Cal
C \subset \Cal F_5^{\text{(big)}} \subset \Cal A_5$. (One uses
the Nash construction (2.1.4) and the fact that the Fermat cubic
appears in the linear system of hyperplane sections of a generic
cubic 3-fold in $\Bbb P^4$.) From the structure of the Gauss map
\cite{C-G, \S\S12-13} one could presumably also infer that $\Cal
C \not\subset \Cal F_5^{\text{(small)}}$. \endremark

\vskip.5\baselineskip \subhead 5.3\quad Questions \endsubhead
\vskip.5\baselineskip

(1) To what extent can the nondegeneracy assumption be removed
from Proposition (2.3.4) on translation hypersurfaces?  In our
work on the Gauss map in genus $g$ (cf\. \cite{MSV2}), we can
allow the curve germs $(\Gamma_i, 0)$ in $\Bbb C^g$ to coincide
in subcollections, say $(\Gamma_i, 0)$ appears $m_i$ times, and
still conclude that $(A, \Theta, \xi) \in \Cal F_g$. But we
assume that the Gauss images $(\overline{\Gamma}_i, p_i)$ in
$\Bbb P^{g-1}$ are smooth disjoint curve germs and the divisor $D
= \sum m_ip_i$ spans a hyperplane of $\Bbb P^{g-1}$. For example,
this situation arises for the theta divisor of the Jacobian of a
smooth genus $4$ canonical curve in $\Bbb P^3$ that has a plane
section of the form $2p + 4q$.  Then $D = p + 2q$, and the cubic
form $\overline{m}_{\xi}(x, y, z)$ attached to the odd theta
characteristic $\xi = \Cal O(D)$ looks like $x^3$, a degenerate
Fermat cubic.

(2) What is the codimension 1 part of the boundary of the closure
of $\Cal F_g$ in a toroidal compactification of $\tilde{\Cal
A}_g$?

(3) Let $\triangle$ be the level $1$ modular form $\prod_{\xi
\text{ odd}} h_{\xi}(\delta)$ and let $D \subset \Cal A_g$ be the
divisor defined by $\triangle = 0$.  Is $D$ irreducible?

(4) Is $\script J_g = \Cal F_g^{(small)} - \Cal R_g$?

\enddocument